\documentclass[reqno,10pt,twoside]{article}
\usepackage[hmargin=1in, vmargin=1.25in, a4paper, centering]{geometry}
\usepackage{amsmath, amsthm, amssymb, enumerate}
\usepackage{fourier}
\usepackage[lining]{ebgaramond}
\usepackage{graphicx}
\usepackage{xcolor}
\usepackage[font=footnotesize, labelfont=bf, labelsep=colon, justification=centerlast, margin=1in]{caption}
\usepackage{etoolbox}
\usepackage[section]{placeins}
\usepackage{titletoc}
\numberwithin{equation}{section}
\theoremstyle{plain}
\newtheorem{theorem}{Theorem}[section]
\newtheorem{catheorem}[theorem]{Computer-assisted Theorem}
\newtheorem{proposition}[theorem]{Proposition}
\newtheorem{lemma}[theorem]{Lemma}
\newtheorem{corollary}[theorem]{Corollary}
\theoremstyle{definition}
\newtheorem{definition}[theorem]{Definition}

\theoremstyle{remark}
\newtheorem{remark}[theorem]{Remark}
\usepackage{etoolbox}
\newcommand{\addQEDstyle}[2]{\AtBeginEnvironment{#1}{\pushQED{\qed}\renewcommand{\qedsymbol}{#2}}\AtEndEnvironment{#1}{\popQED}}
\addQEDstyle{remark}{$\blacktriangleleft$}
\addQEDstyle{example}{$\blacktriangleleft$}
\addQEDstyle{definition}{$\blacktriangleleft$}
\newcommand{\area}{\operatorname{Area}}
\newcommand{\Dir}{{\mathrm{Dir}}}
\newcommand{\dr}{\mathrm{d}}
\newcommand{\ir}{\mathrm{i}}
\newcommand{\er}{\mathrm{e}}

\renewcommand{\epsilon}{\varepsilon}

\newcommand{\Lfloor}{\left\lfloor}
\newcommand{\Rfloor}{\right\rfloor}
\newcommand{\entire}[1]{\Lfloor #1 \Rfloor}

\renewcommand{\phi}{\varphi}
\newcommand{\Lip}{\mathrm{Lip}}
\newcommand{\RL}{\mathrm{R}\Lambda}
\newcommand{\LM}{\Lambda\mathrm{M}}
\renewcommand{\tilde}{\widetilde}
\renewcommand\footnotemark{}
\usepackage{titlesec}
\titleformat{\section}
{\normalfont\large\bfseries}
{\filcenter\S\thesection.}{1ex}{\filcenter}
\titleformat{\subsection}
{\normalfont\bfseries}
{\filcenter\S\thesubsection.}{1ex}{\filcenter}
\usepackage[colorlinks,allcolors=blue,pagebackref]{hyperref}
\usepackage{pifont}
\renewcommand*{\backrefalt}[4]{%
\ifcase #1 %
No citations%
\or
\ding{43}~p.~#2%
\else
\ding{43}~pp.~#2%
\fi}
\newcommand{\mydoi}[1]{\href{https://doi.org/#1}{doi: #1}}
\newcommand{\myarXiv}[1]{\href{https://arxiv.org/abs/#1}{arXiv: #1}}
\begin{document}
\title{%
\vspace{-2cm}P\'{o}lya's conjecture for Dirichlet eigenvalues of annuli
\footnote{{\bf MSC(2020): }Primary 35P15. Secondary 35P20, 33C10, 11P21, 65G20}%
\footnote{{\bf Keywords: }Laplacian, eigenvalues, Weyl's law, lattice points, Bessel functions, Bessel phase functions}%
}
\author{
Nikolay Filonov
\thanks{%
\textbf{N. F.: }St.~Petersburg Department
of Steklov Institute of Mathematics of RAS,
Fontanka 27, 191023, St. Petersburg, Russia;
\href{mailto:filonov@pdmi.ras.ru}{\nolinkurl{filonov@pdmi.ras.ru}}; ORCID: 0000-0002-1586-3031%
}
\and
Michael Levitin\hspace{-3ex}
\thanks{%
\textbf{M. L.: }Department of Mathematics and Statistics, University of Reading, 
Pepper Lane, Whiteknights, Reading RG6 6AX, UK;
\href{mailto:M.Levitin@reading.ac.uk}{\nolinkurl{M.Levitin@reading.ac.uk}}; \url{https://www.michaellevitin.net}; ORCID: 0000-0003-0020-3265%
}
\and 
Iosif Polterovich
\thanks{%
\textbf{I. P.: }D\'e\-par\-te\-ment de math\'ematiques et de statistique, Univer\-sit\'e de Mont\-r\'eal, 
CP 6128 succ Centre-Ville, Mont\-r\'eal QC  H3C 3J7, Canada;
\href{mailto:iosif.polterovich@umontreal.ca}{\nolinkurl{iosif.polterovich@umontreal.ca}}; \url{https://www.dms.umontreal.ca/\~iossif}; ORCID: 0009-0007-0052-6589%
}
\and
David A. Sher
\thanks{%
\textbf{D. A. S.:  }Department of Mathematical Sciences, DePaul University, 2320 N. Kenmore Ave, 60614, Chicago, IL, USA;
\href{mailto:dsher@depaul.edu}{\nolinkurl{dsher@depaul.edu}}; ORCID: 0009-0003-2478-1083%
}
}
\date{\small published in Journal of the London Mathematical Society \textbf{113}:2, e70425 (2026). \mydoi{10.1112/jlms.70425}.} 
\maketitle

\begin{abstract}
We prove P\'olya's conjecture for the eigenvalues of the Dirichlet Laplacian on annular domains. Our approach builds upon and extends the methods we previously developed for disks and balls.  It combines variational bounds, estimates of Bessel phase functions,  refined lattice point counting techniques, and a rigorous computer-assisted analysis. As a by-product, we also derive a two-term upper bound for the Dirichlet eigenvalue counting function of the disk, improving upon P\'olya's original estimate.
\end{abstract}
{\small\tableofcontents}
\section{Introduction and main results}\label{sec:intro}
\subsection{Preliminaries}
Let $\Omega\subset\mathbb{R}^2$ be a bounded  domain,  and let 
\[
\lambda_1(\Omega)\le \lambda_2(\Omega)\le \dots
\]
be the eigenvalues of the Dirichlet Laplacian $-\Delta^\Dir_\Omega$ in $\Omega$, enumerated with multiplicities.  Denote by 
\[
\mathcal{N}^\Dir_\Omega(\lambda):=\#\left\{n\in\mathbb{N}: \lambda_n(\Omega)\le \lambda^2\right\}
\]
the corresponding eigenvalue counting function. By Weyl's law, its  asymptotics  is given by 
\begin{equation}\label{eq:Weyl}
\mathcal{N}^\Dir_\Omega(\lambda)=\frac{\area(\Omega)}{4\pi} \lambda^2+o\left(\lambda^2\right)\qquad\text{as }\lambda\to +\infty.
\end{equation}

More than seventy years ago,  George P\'olya \cite{Polya} conjectured that the leading term of the asymptotics  \eqref{eq:Weyl} is in fact a uniform bound on $\mathcal{N}^\Dir_\Omega(\lambda)$ for any $\Omega\subset \mathbb{R}^2$: namely, that for all $\lambda>0$,
\begin{equation}\label{eq:Polya}
\mathcal{N}^\Dir_\Omega(\lambda) < \frac{\area(\Omega)}{4\pi} \lambda^2.
\end{equation}
P\'olya himself later proved this for domains which \emph{tile} the  plane \cite{Polya2}. No other planar domains for which  P\'olya's conjecture holds were known until very recently,  when the authors of this paper proved it for disks and finite sectors,  see \cite{FLPS}. That paper also gives a more comprehensive account of P\'olya's conjecture for the Dirichlet Laplacian and its counterpart in the Neumann case, as well as its higher-dimensional generalisations. 

\subsection{Main result}
The goal of this article is to extend the proof of P\'olya's conjecture to annuli
\[
A_r=\left\{(x_1,x_2)\in\mathbb{R}^2: r^2<x_1^2+x_2^2<1\right\},
\] 
with inner radius $r\in(0,1)$ and outer radius one. Note that any annulus is of this form up to rescaling. Note also that $r=0$ corresponds to the punctured disk which has the same spectrum as the disk (see also \cite[Example 2.2.23]{LMP}) and, therefore, is covered by the results of  \cite{FLPS}.

For brevity, we will denote the eigenvalues of an annulus and its eigenvalue counting function by
\[
\lambda_{r, n}:=\lambda_n(A_r),\qquad \mathcal{N}_r(\lambda):=\mathcal{N}^\Dir_{A_r}(\lambda).
\]
Given that $\area(A_r)=\pi(1-r^2)$, 
P\'olya's conjecture \eqref{eq:Polya} for annuli  states that the inequality
\begin{equation}\label{eq:PolyaA}
\mathcal{N}_r(\lambda) < \frac{1-r^2}{4} \lambda^2
\end{equation}
holds for all  pairs $(r,\lambda)$ in the parameter space
\[
\RL:=\left\{(r,\lambda): 0<r<1, \lambda>0\right\}.
\]
It will be often convenient to work with the parameter 
\[
\mu:=r\lambda\in[0,\lambda)
\]
instead of $r$. In terms of this parameter, the equivalent form of P\'olya's conjecture for the Dirichlet Laplacian on annuli is that the inequality
\begin{equation}\label{eq:PolyaB}
\mathcal{N}_{\frac{\mu}{\lambda}}(\lambda) < \frac{\lambda^2-\mu^2}{4}
\end{equation}
holds for all pairs $(\lambda,\mu)$  in the parameter space
\[
\LM:=\left\{(\lambda,\mu): 0<\mu<\lambda\right\}.
\]

Our main result is

\begin{theorem}\label{thm:main}
P\'olya's conjecture for the Dirichlet Laplacian holds for all planar annuli. In other words,
inequality \eqref{eq:PolyaA} holds for all $(r,\lambda)\in\RL$, or, equivalently,  \eqref{eq:PolyaB} holds for all $(\lambda,\mu)\in\LM$.
\end{theorem}

Note that Theorem \ref{thm:main} gives a first example of a non-simply connected planar domain for which P\'olya's conjecture has been verified.

\begin{remark} The eigenvalues of planar annuli have also been investigated in \cite{GMWW} and more recently in \cite{Guo}. 
Using number-theoretic techniques, these papers improve remainder estimates in the  two-term Weyl asymptotics for annuli.  
Note that the two-term asymptotics implies that P\'olya's conjecture holds for sufficiently large (but unspecified) 
values of $\lambda$.  While the present  paper shares  with \cite{GMWW, Guo} the method of relating the eigenvalue count to a lattice point count,  our results are in a way complementary, as we focus on accurate eigenvalue bounds rather than asymptotics.  
\end{remark} 

\subsection{Ideas of the proof} \label{sec:ideas}
\begin{figure}[ht]
\centering
\includegraphics{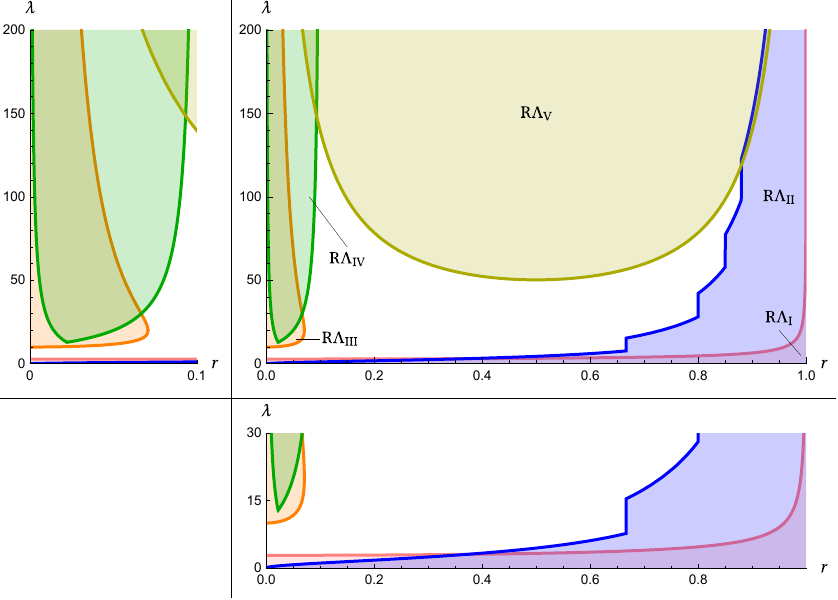}
\caption{Summary of all regions in $(r,\lambda)$-plane for which P\'olya's conjecture is analytically proven to be true, with two side figures zooming near the axes. The exact definitions of the regions can be found in Theorem \ref{thm:regI} for $\RL_\mathrm{I}$, Theorem \ref{thm:vari} for $\RL_\mathrm{II}$, Theorem \ref{thm:III} for $\RL_\mathrm{III}$, Theorem \ref{thm:IV} for $\RL_\mathrm{IV}$, and Theorem \ref{thm:V} for $\RL_\mathrm{V}$. \label{fig:grid}}
\end{figure}

The proof of Theorem \ref{thm:main}  combines  analytic arguments with  a rigorous computer-assisted verification. We analytically prove Theorem \ref{thm:main} in five distinct but partially intersecting regions $\RL_\aleph$, $\aleph=\mathrm{I}, \dots, \mathrm{V}$,  of $\RL$. We postpone their exact definitions until later,  but show their graphical representation in Figure \ref{fig:grid}.

Different approaches are used to deal with each of these regions.  In the  region  
$\RL_\mathrm{I}$ the result follows from the Faber--Krahn and Krahn--Szego isoperimetric inequalities for the first and the second Dirichlet eigenvalues.  In the  region $\RL_\mathrm{II}$,  we use a comparison between the counting functions for the  annuli and for certain carefully chosen flat cylinders, see Lemma \ref{lem:NN}. We then take advantage of the fact that the eigenvalues of these cylinders can be computed explicitly, see  Theorem \ref{thm:vari}.  Note that the arguments in the first two regions are of purely spectral-theoretic nature.

The proofs in the remaining three regions are based on  results concerning the Bessel  phase functions as well as  lattice-counting techniques. Here we extend the approach developed in \cite{FLPS} for the disk. In order to adapt it to the case of annuli we have to overcome several new technical challenges. 

First, instead of  estimating the Bessel phase function itself, we now need to bound the difference between the values of the phase functions at two distinct points, which are chosen depending on the spectral parameter and the inner radius of the annulus (recall that the outer radius is fixed and is equal to one).  Some of these bounds follow from the results of  \cite{FLPSbessel}, but some are new, see  \S\ref{sec:gammabounds}.

Second, we have to use  a more elaborate reduction to a lattice counting problem. In the region $\RL_\mathrm{III}$ we count points below a single curve as in \cite{FLPS}.  In the other two regions, in order to get accurate enough bounds,  we need to split the interval $[0, \lambda]$ into two (in $\RL_\mathrm{IV}$) and three (in  $\RL_\mathrm{V}$) subintervals, respectively, and in each subinterval we count integer points under different  curves, depending on which bound for the phase functions we are using (see Lemma \ref{cor:deltabound}
and Figure \ref{fig:boundsGFH}). 

Third, we need to improve upon  the lattice counting techniques of \cite{FLPS}. Recall that 
\cite[Theorem 5.1]{FLPS} provided  a way to estimate  the number of lattice points under a graph of a non-negative decreasing convex function with a Lipschitz constant $\frac{1}{2}$. We need to extend this bound to concave functions, as well as  to relax the assumption on the Lipschitz constant.  We obtain several  results in this direction in \S\ref{sec:sums}. We believe that they could be of independent interest, and for illustration we present one of them below. 

The following definition is central for our analysis.
\begin{definition}\label{def:tfs}
Let $a,b\in\mathbb{Z}$ with $a<b$, and let $g:[a,b]\to\mathbb{R}$. The \emph{trapezoidal floor sum} of $g$ on $[a,b]$ is defined as
\[
\mathbf{T}(g,a,b) = \frac{1}{2} \sum_{m=a}^{b-1}\left(\entire{g(m)} + \entire{g(m+1)}\right)
= \frac{1}{2}\entire{g(a)}+\sum_{m=a+1}^{b-1}\entire{g(m)} + \frac{1}{2}\entire{g(b)}.
\]
\end{definition}

For a  non-negative function $g$, the trapezoidal floor sum $\mathbf{T}(g,a,b)$ counts the number of  points of the lattice $\mathbb{Z}\times\mathbb{N}$ lying on or under the graph of $g$ on the interval $[a,b]$, with the exception that the lattice points with the edge abscissas $a$ and $b$ are counted with the weight $\frac{1}{2}$. A useful observation is that for an integer $p \in(a,b)$, 
\[
\mathbf{T}(g,a,b)=\mathbf{T}(g,a,p)+\mathbf{T}(g,p,b).
\]
For example, we prove the following estimate on the trapezoidal floor sum for concave Lipschitz functions that is used in the region $\RL_\mathrm{V}$.

\begin{theorem}\label{thm:t25intro}
Let $\alpha, \beta \in \mathbb{Z}$, $\alpha<\beta$.  Let $g$ be a decreasing concave function on $[\alpha,\beta]$ 
which is Lipschitz with a constant $c \in(0,1)$, and assume additionally that
\[
\entire{g(\alpha)}=\entire{g(p)} > \entire{g(p+1)}
\]
for some integer $p \in[\alpha,\beta)$.
Then
\[
\mathbf{T}(g,\alpha,\beta) \le \int_\alpha^\beta g(z)\,\dr z - \frac{1-c}{2}(\beta-p).
\]
\end{theorem}

Note that without the second term on the right,  the bound is an easy consequence of concavity only, see Proposition \ref{prop:gconcave}. The proof of Theorem \ref{thm:t25intro} can be found in \S\ref{sec:sums}. For some related asymptotic estimates see also \cite{LL1, LL2}.

For decreasing convex functions we prove a related result under the assumption that the function is Lipschitz with a constant $\frac{1}{2}$ on the whole interval, and {\em additionally} is Lipschitz with a constant $\frac{1}{3}$ on a subinterval, see Theorem \ref{thm:conveximproved}.  One consequence of this is the following improvement of P\'olya's conjecture for the Dirichlet Laplacian on the disk \cite{FLPS}.

\begin{theorem}\label{thm:diskimproved}
Let $\mathbb{D}$ be a unit disk. Then the following two-term P\'olya-type bound for the Dirichlet Laplacian holds,
\begin{equation}\label{eq:improvedPolya}
\mathcal{N}^\Dir_\mathbb{D}(\lambda) < \frac{\lambda^2}{4}-\frac{\entire{\omega_0\lambda}}{2},
\end{equation} 
where 
\[
\omega_0:=\frac{\sqrt{3}}{2\pi}-\frac{1}{6}\approx 0.108998.
\]
\end{theorem}

We prove Theorem \ref{thm:diskimproved} in \S\ref{subsec:regionIII}.

\begin{remark} While the present paper was under review,  estimate \eqref{eq:improvedPolya}  was further improved in \cite[Theorem 1.7]{Guo-p} to
\[
\mathcal{N}^\Dir_\mathbb{D}(\lambda) < \frac{\lambda^2}{4}-\frac{\entire{\omega_0\lambda}}{8}-\frac{3\entire{\omega_1\lambda}}{8},
\]
where
\[
\omega_1:=\frac{\sqrt{2+\sqrt{2}}}{2\pi} - \frac{3\sqrt{2-\sqrt{2}}}{16}\approx 0.150574.
\]
\end{remark}

\begin{figure}[ht]
\centering
\includegraphics{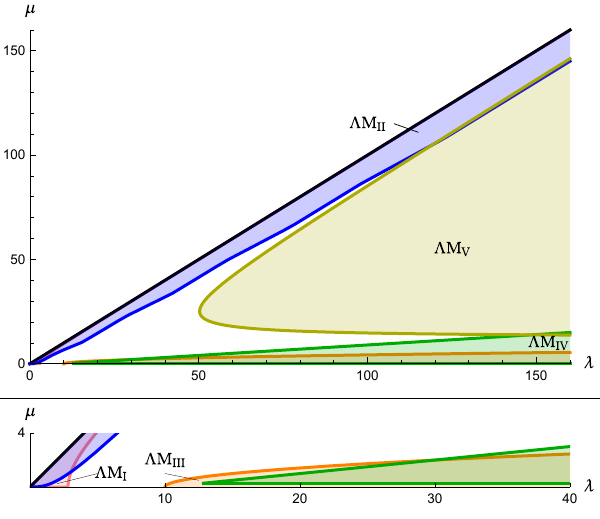}
\caption{Summary of all regions in $(\lambda, \mu)$-plane for which P\'olya's conjecture is analytically proven to be true. The bottom figure zooms near the $\lambda$-axis and the origin. \label{fig:LM}}
\end{figure}

The analytic results establishing  P\'olya's conjecture for the regions $\RL_\aleph$ can be equivalently formulated  for the regions 
$\LM_\aleph$, $\aleph=\mathrm{I}, \dots \mathrm{V}$, see Figure \ref{fig:LM}. 
In fact, such a reformulation is more convenient  for the computer-assisted part of  the proof of Theorem \ref{thm:main}. 

Set 
\[
\LM_\mathrm{theory} := \bigcup_{\aleph =\mathrm{I}}^\mathrm{V} \LM_\aleph.
\]
Importantly, the remaining region 
\[
\LM\setminus\LM_\mathrm{theory}
\]
is bounded, see Theorem \ref{thm:finiteQ}, and we implement  a rigorous computer-assisted algorithm verifying P\'olya's conjecture there. Our method builds upon the approach of \cite[\S8]{FLPS}. We show that if a required inequality on the lattice count holds at a given point with some margin, then it holds in a certain explicitly described domain around this point, see Lemma \ref{lem:rect}. This observation, together with the technique of verified rational approximations,  allows us to rigorously check P\'olya's conjecture in a region $\LM_{\mathrm{comp}}\supset \LM\setminus~\LM_\mathrm{theory}$. 
This completes the proof of Theorem \ref{thm:main}.

\subsection{Plan of the paper} The paper is organised as follows. In \S\ref{sec:IandII},  we prove P\'olya's conjecture in the region $\RL_\mathrm{I}$ using isoperimetric inequalities for the first two Dirichlet eigenvalues,  and  in the region $\RL_\mathrm{II}$ via comparison with the spectra of flat cylinders.
In \S\ref{sec:sums}, we present abstract results on the trapezoidal floor sums for concave and convex functions satisfying different assumptions on their Lipschitz constants. In \S\ref{sec:eigBessel}, we relate the eigenvalue counting function of an annulus to the difference of the Bessel phase functions. Bounds on these differences in terms of some elementary functions and  a reduction to a lattice counting problem are presented in \S\ref{sec:gammabounds}. Some useful auxiliary properties of these elementary functions are collected in \S\ref{sec:aux}. In \S\ref{sec:apptrap}, we apply the results on the trapezoidal floor sums and the bounds on the difference of the Bessel phase functions in order to verify 
P\'olya's conjecture in the regions $\RL_\mathrm{III}, \RL_\mathrm{IV}$ and $\RL_\mathrm{V}$. This section is the most technical part of the paper. Finally, in \S\ref{sec:compass} we present a rigorous computer-assisted algorithm which affirms P\'olya's conjecture in the remaining part of $\RL$.

\section{Regions I and II via  isoperimetric inequalities and comparison with flat cylinders}\label{sec:IandII} 

By the Faber--Krahn and Krahn--Szego inequalities (see, e.g., \cite{Henrot, LMP}), we know that P\'olya's conjecture is  true for any domain $\Omega\subset\mathbb{R}^d$  for values of $\lambda$ such that  the corresponding Weyl's term satisfies $C_d|\Omega|\lambda^2\le 2$, which in the case of annulus $A_r$ is equivalent to
\[
\lambda^2-8\le \mu^2.
\]
Therefore, we arrive at

\begin{theorem}\label{thm:regI} 
Let
\[
\eta_\mathrm{I}(r) := \sqrt{\frac{8}{1-r^2}},\qquad \zeta_\mathrm{I}(\lambda):=\sqrt{\lambda^2-8}.
\]
Inequality \eqref{eq:PolyaB}  holds for all $(\lambda,\mu)\in\LM_\mathrm{I}$, where
\[
\LM_\mathrm{I}:=\left\{(\lambda,\mu): 0<\mu<\lambda\le 2\sqrt{2}\right\}\cup\left\{(\lambda,\mu): \lambda>2\sqrt{2},\ \zeta_\mathrm{I}(\lambda)\le \mu<\lambda\right\}\subset\LM.
\]
Equivalently,  inequality \eqref{eq:PolyaA}  holds for  all $(r,\lambda)\in\RL_\mathrm{I}$, where
\[
\RL_\mathrm{I}:=\left\{(r,\lambda): 0<r<1,\ \lambda\le \eta_\mathrm{I}(r)\right\}\subset\RL.
\]
\end{theorem}
 
We can extend the result of Theorem \ref{thm:regI} by a more involved argument. In order to do so, consider first the Dirichlet Laplacian on a flat cylinder $\mathcal{C}_h:=(1-h,1)\times \mathbb{S}^1$, with $h>0$ (cf. \cite{FS} where some related results were obtained). An elementary separation of variables shows that the (unordered) eigenvalues of $-\Delta_{\mathcal{C}_h}^\Dir$ are given by 
\[
\widetilde{\lambda}_{n,m}\left(\mathcal{C}_h\right)=m^2+\frac{\pi^2 n^2}{h^2},\qquad  (n,m)\in\mathbb{N}\times\mathbb{Z}.
\] 
The corresponding eigenvalue counting function is
\[
\widetilde{\mathcal{N}}_h(\lambda):=\mathcal{N}_{\mathcal{C}_h}(\lambda)=\#\left\{(n,m)\in\mathbb{N}\times\mathbb{Z}: m^2+\frac{\pi^2 n^2}{h^2}\le \lambda^2\right\}.
\]

Recall that the eigenvalue counting function of the annulus $A_r$ is denoted by $\mathcal{N}_r(\lambda)$. We have

\begin{lemma}\label{lem:NN} Let $0<r<1$, and set $h=h_r:=\frac{1-r}{\sqrt{r}}$. 
Then 
\[
\mathcal{N}_r(\lambda)\le \widetilde{\mathcal{N}}_h(\lambda)\qquad\text{for all }\lambda>0.
\]
\end{lemma} 
 
\begin{proof} Let 
\begin{equation}\label{eq:fseries}
f(\rho,\psi)=\sum_{m\in\mathbb{Z}} f_m(\rho)\er^{\ir m\psi}
\end{equation} 
be a non-zero function from  the Sobolev space $H^1_0(A_r)$, so that $f_m(r)=f_m(1)=0$. The Rayleigh quotient of the operator $-\Delta^\Dir_{A_r}$ acting on $f$ is
\begin{equation}\label{eq:RAr}
\mathtt{R}_r[f]=\frac{\sum\limits_{m\in\mathbb{Z}} \int\limits_{r}^1 \left(\rho\left|f_m'(\rho)\right|^2 +\frac{m^2}{\rho}\left|f_m(\rho)\right|^2\right)\,\dr\rho}{\sum\limits_{m\in\mathbb{Z}}\int\limits_{r}^1 \left|f_m(\rho)\right|^2 \rho \,\dr\rho}.
\end{equation}

At the same time, the eigenvalues of $-\Delta^\Dir_{\mathcal{C}_{h_r}}$ coincide with those of the operator 
\begin{equation}\label{eq:scale}
-r\frac{\partial^2}{\partial \rho^2}-\frac{\partial^2}{\partial \psi^2}
\end{equation}
acting in $H^1_0\left(\mathcal{C}_{1-r}\right)$ (with the cylindrical coordinates $(\rho, \psi)$ coinciding with the polar coordinates in $A_r$), see Figure \ref{fig:anncylalt}. 

\begin{figure}[ht]
\centering
\includegraphics{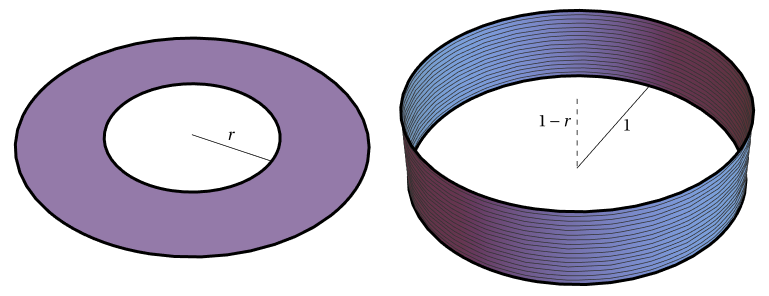}
\caption{An annulus $A_r$ and a cylinder $\mathcal{C}_{1-r}$. Note that we are comparing the spectrum of the Dirichlet Laplacian in $A_r$ with the spectrum of the Dirichlet realisation of \eqref{eq:scale} in $\mathcal{C}_{1-r}$, which in turn coincides with the spectrum of the Dirichlet Laplacian in $\mathcal{C}_{h_r}$. \label{fig:anncylalt}}
\end{figure}

Note that  $f$  is an element of  $H^1_0\left(\mathcal{C}_{1-r}\right)$ if and only if it belongs to $H^1_0\left(A_r\right)$, and the Rayleigh quotient of the operator \eqref{eq:scale} acting on \eqref{eq:fseries} is
\begin{equation}\label{eq:RCh}
\widetilde{\mathtt{R}}_r[f]
=\frac{\sum\limits_{m\in\mathbb{Z}} \int\limits_{r}^1 \left(r \left|f_m'(\rho)\right|^2+ m^2 \left|f_m(\rho)\right|^2\right)\,\dr\rho}{\sum\limits_{m\in\mathbb{Z}}\int\limits_{r}^1 \left|f_m(\rho)\right|^2\,\dr\rho}.
\end{equation}

The comparison of the integrands in \eqref{eq:RCh} and \eqref{eq:RAr} immediately gives $\mathtt{R}[f] > \widetilde{\mathtt{R}}_r[f]$ for any $f\in H^1_0\left(A_r\right)\setminus\{0\}$, 
and the result follows from the variational principle.
\end{proof} 
 
\begin{figure}[ht]
\centering
\includegraphics{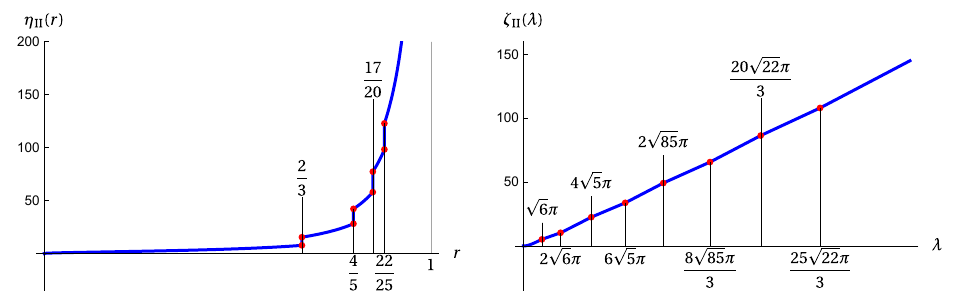}
\caption{The plots of $\eta_\mathrm{II}(r)$ and $\zeta_\mathrm{II}(\lambda)$. The red dots indicate the positions of singularities in definitions \eqref{eq:eta2} and \eqref{eq:zeta2}.\label{fig:etazetaII}}
\end{figure}

\begin{theorem}\label{thm:vari}
Let
\begin{equation}\label{eq:rjs}
r_0:=0, \quad r_1:=\frac{2}{3},\quad r_2:=\frac{4}{5}, \quad r_3:=\frac{17}{20}, \quad r_4:=\frac{22}{25},\quad  r_5:=1,
\end{equation}
and set
\begin{equation}\label{eq:eta2}
\eta_\mathrm{II}(r):= \frac{(j+1)\pi\sqrt{r}}{1-r}\quad\text{ if}\quad r_{j}\le r< r_{j+1},\qquad j=0,1, 2, 3, 4,
\end{equation}
\begin{equation}\label{eq:zeta2}
\zeta_\mathrm{II}(\lambda):=\begin{cases} 
\lambda-\frac{i\pi}{2\lambda}\left(\sqrt{4\lambda^2+i^2\pi^2}-i\pi\right)&\text{ if}\quad  i\pi\frac{\sqrt{r_{i-1}}}{1-r_{i-1}}\le \lambda< i\pi\frac{\sqrt{r_i}}{1-r_i},\qquad i=1, 2, 3, 4, 5,\\
 r_i \lambda &\text{ if}\quad i\pi\frac{\sqrt{r_i}}{1-r_i} \le \lambda < (i+1)\pi\frac{\sqrt{r_i}}{1-r_i},\qquad i=1, 2, 3, 4,
\end{cases}
\end{equation}
see Figure \ref{fig:etazetaII}.

Inequality \eqref{eq:PolyaA} holds for all $(r,\lambda)\in\RL_\mathrm{II}$, where
\begin{equation}\label{eq:reg2r}
\RL_\mathrm{II}:=\left\{(r, \lambda): 0<r<1,\ 0<\lambda<\eta_\mathrm{II}(r)
\right\}.
\end{equation}
Equivalently, inequality \eqref{eq:PolyaB}  holds for all $(\lambda, \mu)\in\LM_\mathrm{II}$, where
\begin{equation}\label{eq:reg2}
\LM_\mathrm{II} = \left\{(\lambda,\mu): \zeta_\mathrm{II}(\lambda)<\mu<\lambda\right\}\subset\LM.
\end{equation}
\end{theorem}

\begin{proof} Fix $r\in(0,1)$. With $h=h_r$, we have
\[
\begin{split}
\widetilde{\mathcal{N}}_{h_r}(\lambda) &= \#\left\{(n,m)\in\mathbb{N}\times\mathbb{Z}: n^2+\frac{(1-r)^2}{r\pi^2}m^2\le \frac{(1-r)^2}{r\pi^2}\lambda^2\right\}\\
&=
\sum_{n=1}^{\entire{\frac{1-r}{\pi\sqrt{r}}\lambda}} \left(1+2\entire{\sqrt{\lambda^2-\frac{\pi^2 r}{(1-r)^2}n^2}}\right),
\end{split}
\]
where we assumed the notational convention $\sum\limits_{n=1}^0 := 0$. By Lemma \ref{lem:NN}, inequality \eqref{eq:PolyaA} will hold for a given $r$ and $\lambda$ if we can show that 
\begin{equation}\label{eq:PolyatildeN}
\widetilde{\mathcal{N}}_{h_r}(\lambda)<\frac{1-r^2}{4}\lambda^2.
\end{equation}

Assume, first of all, that $\entire{\frac{1-r}{\pi\sqrt{r}}\lambda}=0$, that is, 
\[
\lambda< \frac{\sqrt{r}\pi}{1-r}.
\]
Then $\widetilde{\mathcal{N}}_{h_r}(\lambda)=0$, and \eqref{eq:PolyatildeN} follows immediately.

Let now $\entire{\frac{1-r}{\pi\sqrt{r}}\lambda}=j\in \mathbb{N}$, that is,
\begin{equation}\label{eq:lambdaj}
j\frac{\sqrt{r}\pi}{1-r}\le \lambda< (j+1)\frac{\sqrt{r}\pi}{1-r}. 
\end{equation} 
In practice, we will only consider cases when $j\le 4$, but the general scheme may be extended further.
Let 
\[
L_n =  \entire{\sqrt{\lambda^2-\frac{\pi^2 r}{(1-r)^2}n^2}},\qquad n=1,\dots,j,
\]
so that 
\begin{equation}\label{eq:Lnineq}
L_n^2 + \frac{\pi^2 r}{(1-r)^2}n^2 \le \lambda^2 < (L_n+1)^2 + \frac{\pi^2 r}{(1-r)^2}n^2.
\end{equation}
Then 
\[
\widetilde{\mathcal{N}}_{h_r}(\lambda) = j +2 L_1+\dots+ 2 L_j.
\]

In order to effectively obtain some restrictions on $r$ which guarantee that inequality \eqref{eq:PolyatildeN} holds we will replace it by a stronger inequality based on the lower bounds on $\lambda^2$ from \eqref{eq:Lnineq}:
\begin{equation}\label{eq:Lnsum}
\begin{split}
\frac{1-r^2}{4}\lambda^2 - \widetilde{\mathcal{N}}_{h_r}(\lambda) 
&=\sum_{n=1}^j \left(\tau_n \frac{1-r^2}{4}  \lambda^2 - 1 - 2 L_n\right) \\
&\ge \sum_{n=1}^j \left(\tau_n \frac{1-r^2}{4} \left(L_n^2 + \frac{\pi^2 r}{(1-r)^2}n^2\right) - 1 - 2 L_n\right) >0,
\end{split}
\end{equation}
where
\begin{equation}\label{eq:taus}
(\tau_1,\dots,\tau_j)=\boldsymbol{\tau}\in(0,1]^j\qquad\text{with } \tau_1+\dots+\tau_j=1,
\end{equation}
are some constants to be chosen later.  Each summand in the right-hand side of \eqref{eq:Lnsum} is minimised, over $L_n\in\mathbb{R}$, by taking $L_n=\frac{4}{(1-r^2)\tau_n}$.  Substituting these values into \eqref{eq:Lnsum} and multiplying the resulting inequality by $4(1-r^2)$, we obtain the inequality 
\begin{equation}\label{eq:psi}
S_j(r; \boldsymbol{\tau}):=\left(\sum_{n=1}^j n^2\tau_n\right)\pi^2 r (1+r)^2 -16\left(\sum_{n=1}^j  \frac{1}{\tau_n}\right) - 4 j (1-r^2)>0.  
\end{equation}
For every fixed choice of a vector $\boldsymbol{\tau}$ satisfying constraints \eqref{eq:taus}, the cubic (in $r$) polynomial $S_j(r; \boldsymbol{\tau})$ is monotone increasing for $r\ge 0$ with $S_j(0; \boldsymbol{\tau})<0$.  
Therefore, it  has a single positive root $r_j^*(\boldsymbol{\tau})$. If we can choose a vector $\boldsymbol{\tau}^+$ and a number $r^+\in(0,1)$ such that 
\[
S_j\left(r^+; \boldsymbol{\tau}^+\right)>0,
\]
then $r^+>r_j^*(\boldsymbol{\tau}^+)$, and so inequality \eqref{eq:psi} with $\boldsymbol{\tau}=\boldsymbol{\tau}^+$, and therefore also inequality \eqref{eq:PolyatildeN},  holds for $r\ge r^+$ 
and for $\lambda$ satisfying \eqref{eq:lambdaj}. 

We now finish the proof of the theorem by looking at the specific cases $j\in\{1, 2, 3, 4\}$. In each case, we choose $r^+:=r_j$, which is given by \eqref{eq:rjs}.
\begin{description}
\item[Case $j=1$.] The only possible choice is $\boldsymbol{\tau}^+=(1)$, and we  get
\[
S_1(r_1; \boldsymbol{\tau}^+)=\frac{2}{27}\left(25\pi^2-246\right)\approx 0.0548>0.
\]

\item[Case $j=2$.]  We choose $\boldsymbol{\tau}^+=\left(\frac{3}{8}, \frac{5}{8}\right)$, cf.\ Figure \ref{fig:psi23} (left), yielding 
\[
S_2(r_2; \boldsymbol{\tau}^+)=\frac{23}{750} \left(243 \pi ^2-2320\right)\approx 2.4016>0.
\]

\begin{figure}[ht]
\centering
\includegraphics{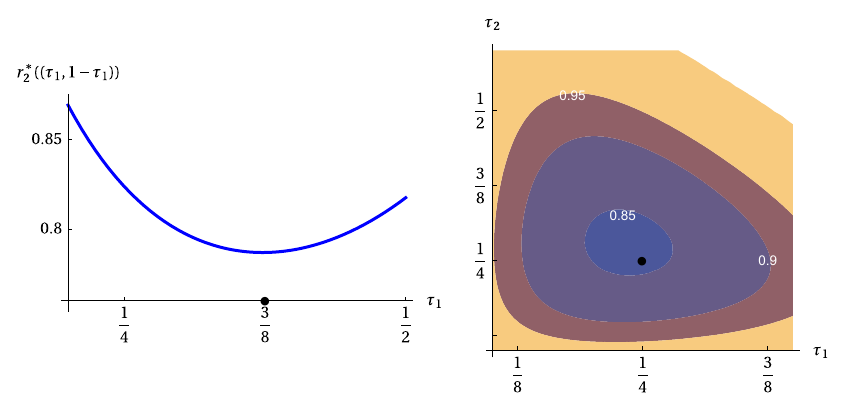}
\caption{On the left, the plot of numerically computed $r_2^*\left((\tau_1,1-\tau_1)\right)$ as a function of $\tau_1$. On the right, the numerically computed contour plot of $r_3^*\left((\tau_1,\tau_2,1-\tau_1-\tau_2)\right)$ in $\left(\tau_1, \tau_2\right)$-coordinates. \label{fig:psi23}}
\end{figure}

\item[Case $j=3$.]  We choose $\boldsymbol{\tau}^+=\left(\frac{1}{4},\frac{1}{4},\frac{1}{2}\right)$, cf.\ Figure \ref{fig:psi23} (right), as 
\[
S_3(r_3; \boldsymbol{\tau}^+)=\frac{535279 \pi ^2-5226560}{32000}\approx 1.7635>0.
\]

\item[Case $j=4$.] We choose, on the basis of some numerical experiments,  $\boldsymbol{\tau}^+=\left(\frac{1}{6}, \frac{1}{6}, \frac{1}{5}, \frac{7}{15}\right)$.
We get 
\[
S_4(r_4; \boldsymbol{\tau}^+)=\frac{17179393\pi^2-169474000}{546875}\approx 0.1459>0.
\] 
\end{description}

The combination of all the  cases completes the proof of \eqref{eq:reg2r}. It remains to remark that  \eqref{eq:reg2} is just an explicit representation of $\left\{(\lambda,\mu): \left(\frac{\mu}{\lambda},\lambda\right)\in \RL_\mathrm{II}\right\}$, where we have taken extra care of the jumps of $\eta_\mathrm{II}(r)$ at $r=r_j$, $j=1, 2, 3, 4$.
\end{proof}

\section{Trapezoidal floor sums for Lipschitz functions and their bounds}\label{sec:sums}

We want to compare, under various conditions on a function $g$, its trapezoidal floor sum with the integrals $\int_a^b g(x)\,\dr x$. Some estimates of this type have been already established in \cite{FLPS} and \cite{FLPS-AB}, however proving Theorem \ref{thm:main} requires more delicate upper bounds for trapezoidal floor sums  which we collect here.

We recall that a function $g:[a,b]\to\mathbb{R}$ is called \emph{$c$-Lipschitz} ($c>0$) if $|g(z)-g(w)|\le c|z-w|$ for all $z,w\in[a,b]$. We will write this as $g\in\Lip_c$.

The following introductory result is obvious.

\begin{proposition}\label{prop:gconcave} If $g(z)$ is concave, then we always have
\[
\mathbf{T}(g,a,b)\le\int_a^b g(z)\,\dr z.
\]
\end{proposition}

\begin{proof}
In this case, $\mathbf{T}(g,a,b)$ is trivially less than or equal to the trapezoid rule approximation for the integral of $g$, which is a lower bound for the true integral. 
\end{proof}

A sharper bound for a concave function $g$ can be obtained subject to some additional constraints. 

\begin{theorem}\label{thm:t25}
Let $a, b \in \mathbb{Z}$, $a<b$. Let $c\in(0,1)$.
Assume that $g$ is a decreasing concave $\Lip_c$ function on $[a,b]$ such that 
\[
\entire{g(a)} > \entire{g(a+1)}.
\]
Then
\begin{equation}\label{eq:gconclip}
\mathbf{T}(g,a,b) \le \int_a^b g(z)\,\dr z - \frac{1-c}{2}(b-a).
\end{equation}
\end{theorem}

\begin{proof} By the conditions of the Theorem, the function $g$ is strictly monotone on some interval to the left of $a+1$, and, therefore by concavity, also on $[a+1,b]$, and thus is invertible there. In addition,
\[
N:=\entire{g(a)} >\entire{g(b)}  =: M.
\]
Denote
\[
a _n = \max \left\{m \in \mathbb{Z}: g(m) \ge n\right\}=\entire{g^{-1}(n)}, \qquad n=N-1, \dots, M+1.
\]
Set also $a_N := a$ and $a_M:=b$. 

By additivity of the trapezoidal floor sums, we have
\[
\mathbf{T}(g,a,b) = \mathbf{T}(g, a_N, a_{N-1}) + \dots + T(g, a_{M+1}, a_M),
\]
and it is therefore sufficient to prove the statement of the Theorem with $(a,b)$ replaced by $(a_n, a_{n-1})$, $n=N,\dots,M-1$. 

Consider any such interval, and set $A:=a_n$, $B:=a_{n-1}$, and $k:=B-A\in\mathbb{N}$. Without loss of generality we can assume $n=1$, otherwise we just need to add the same constant $(n-1)k$ to both sides of \eqref{eq:gconclip} which does not affect the validity of the inequality. Then $g(A)\ge 1>g(A+1)>\dots>g(B)\ge 0$ and $\mathbf{T}(g, A, B)=\frac{1}{2}$. We therefore need to prove that 
\begin{equation}\label{5*}
\int_A^B g(z)\,\dr z \ge \frac{1}{2} + \frac{k(1-c)}{2} = \frac{k(1-c)+1}{2}.
\end{equation}
By concavity of $g$ and using $g(A)\ge 1$, $g(B)\ge 0$, we have
\[
\int_A^B g(z)\,\dr z - \frac{k(1-c)+1}{2} \ge \frac{(g(A)+g(B))k}{2} -  \frac{k(1-c)+1}{2} \ge \frac{k}{2} -  \frac{k(1-c)+1}{2} = \frac{ck-1}{2},
\]
and therefore \eqref{5*} holds if $k\ge\frac{1}{c}$. 

Suppose now $k\le \frac{1}{c}$. We have, by the Lipschitz condition on $g$, the bound $g(B)\ge 1 - ck$, and  therefore
\[
\begin{split}
\int_A^B g(z)\,\dr z  - \frac{k(1-c)+1}{2} &\ge \frac{(2 - c k)k}{2} - \frac{k(1-c)+1}{2} = \frac{2 k  -c k^2 + c k - k -1}{2} \\&=  \frac{(k-1)(1-c k)}{2}\ge 0.
\end{split}
\]
Thus, \eqref{5*} holds again.
\end{proof}

We can now easily give
\begin{proof}[Proof of Theorem \ref{thm:t25intro}] 
By additivity, 
\[
\mathbf{T}(g,\alpha,\beta)=\mathbf{T}(g,\alpha,p)+\mathbf{T}(g,p,\beta),
\]
and we estimate the first term by Proposition \ref{prop:gconcave} with $[a,b]=[\alpha,p]$, and the second term by Theorem \ref{thm:t25} with $[a,b]=[p,\beta]$.
\end{proof}

We now switch to the bounds for convex decreasing functions $g$. The next result is essentially a version of \cite[Theorem 5.1]{FLPS}. 

\begin{theorem}\label{thm:convexold} Let $a,b\in\mathbb{Z}$, and let $g:[a,b]\to \mathbb{R}$ be non-negative, decreasing, convex, and of class $\Lip_{\frac{1}{2}}$. Suppose also that $g(b)$ is an integer. Then
\[
\mathbf{T}\left(g+\frac{1}{4},a,b\right)\le\int_a^b g(z)\,\dr z.
\]
Moreover, equality implies that $g(z)$ is a constant on $[a,b]$.
\end{theorem}

\begin{proof} If $g(b)=0$, this is an immediate consequence of \cite[Theorem 5.1]{FLPS}, noting that the last term in the sum is zero and so it does not matter whether it has a factor of $1/2$. Adding an integer $n$ to $g(z)$ increases both sides by $(b-a)n$ and so equality is maintained.
\end{proof}

The bound in Theorem \ref{thm:convexold} can be improved given more information on the function $g$.

\begin{theorem}\label{thm:conveximproved}  Let $a,b\in\mathbb{Z}$, and suppose $g(z)$ is decreasing, convex,  $\Lip_{\frac{1}{2}}$ on $[a,b]$, and has $g(b)=0$. Suppose further that there exists $t\in [a,b]$ for which $g(z)$ is $\Lip_{\frac{1}{3}}$ on $[t,b]$. Then
\[
\mathbf{T}\left(g+\frac{1}{4},a,b\right)\le \int_a^b g(z)\,\dr z - \frac{1}{4}\entire{g(t)}.
\]
\end{theorem}

In order to prove  Theorem \ref{thm:conveximproved} we require

\begin{lemma}\label{lem:improved} Let $A,B\in\mathbb{Z}$. Suppose that $g(z)$ is decreasing and convex on $[A,B+1]$ with Lipschitz constant $\frac{1}{2}$, and is Lipschitz with constant $\frac{1}{3}$ for $z\ge A+1$. Suppose further that for some $n\in\mathbb Z$,
\[
g(A)\ge n+1\ge g(A+1)\ge g(B)\ge n\ge g(B+1).
\]
Then
\begin{equation}\label{eq:lemmanikolay}
\mathbf{T}\left(g+\frac{1}{4}, A, B\right)\le\int_A^B g(t)\, \dr t - \frac{1}{4}.
\end{equation} 
\end{lemma}

\begin{proof}[Proof of Lemma \ref{lem:improved}] Without loss of generality we assume $A=0$ and $n=0$. Let 
\[
k=\#\{m\in[0,B]:g(m)\ge 3/4\},
\] 
and consider three cases.

\textbf{Case 1:} $k=1$. The left-hand side of \eqref{eq:lemmanikolay} is $\frac{1}{2}$. Since $g(0)\ge 1$, we must have $g(z)\ge 1-\frac{z}{2}$ for $z\ge 0$. In particular $g(1)\ge\frac{1}{2}$, and  therefore $g(2)\ge\frac{1}{6}$, and so $B\ge 2$. Thus
\[
\int_0^B g(z)\,\dr z \ge \int_0^2 \left(1-\frac{z}{2}\right)\,\dr z = 1,
\]
and so the right-hand side of \eqref{eq:lemmanikolay} is at least $\frac{3}{4}>\frac{1}{2}$ as desired.

\textbf{Case 2:} $k=2$. The left-hand side of \eqref{eq:lemmanikolay} is $\frac{3}{2}$. Since $g(1)\ge \frac{3}{4}$, we have $g(z)\ge \frac{3}{4}-\frac{z-1}{3}$ whenever $z\ge 1$. Thus, $g(3)\ge\frac{1}{12}$ and  $B\ge 3$. Now the right-hand side of \eqref{eq:lemmanikolay} is at least
\[
\int_0^3 g(z)\,\dr z = \int_0^2 g(z)\,\dr z + \int_2^3 g(z)\,\dr z.
\]
By convexity, the first integral is at least $2g(1)\ge\frac{3}{2}$. The second integral, using the bound $g(z)\ge \frac{3}{4}-\frac{z-1}{3}$, is at least $\frac{1}{4}$. So the right-hand side of \eqref{eq:lemmanikolay} is at least $\frac{3}{2}$ as desired.

\textbf{Case 3:} $k\ge 3$. The left-hand side of \eqref{eq:lemmanikolay} is $k-\frac{1}{2}$. Now we know that $g(k-1)\ge\frac{3}{4}$. By convexity, $g(1)+g(2k-3)\ge 2g(k-1)\ge \frac{3}{2}$. Since $g(1)<1$, we have $g(2k-3)\ge \frac{1}{2}$ and therefore of course $B\ge 2k-2$. So,
\[
\int_0^B g(z)\,\dr z \ge \int_0^{2k-2}g(z)\,\dr z,
\]
which by convexity is at least $(2k-2)g(k-1)\ge \frac{3}{2}(k-1)\ge k$. Thus,  the right-hand side of \eqref{eq:lemmanikolay} is at least $k-\frac{1}{4} > k-\frac{1}{2}$, completing the proof.
\end{proof}

\begin{proof}[Proof of Theorem \ref{thm:conveximproved}]  
Let $N=\entire{g(a)}$, and for each $n\in[0,N]$, let $q_n=\max\{m\in\mathbb Z : g(m)\ge n\}$. Then
\[
\mathbf{T}\left(g+\frac{1}{4},a,b\right)=\mathbf{T}\left(g+\frac{1}{4},a,q_N\right) + \sum_{n=0}^{N-1}\mathbf{T}\left(g+\frac{1}{4},q_{n+1},q_n\right).
\]
Apply Lemma \ref{lem:improved} for each $n$ for which $q_{n+1}+1\ge t$, with $A=q_{n+1}$ and $B=q_n$, and use Theorem \ref{thm:convexold} to estimate each of the other sums. We get the integral we want plus an additional $-\frac{1}{4}$ for every $n$ for which $q_{n+1}+1\ge t$. Since $g(t) \le g(q_{n+1}+1) < n+1$ for any such $n$, there are at least $\entire{g(t)}$ of these values of $n$, which completes the proof.
\end{proof}

\section{The eigenvalue counting function and the Bessel phase function}
\label{sec:eigBessel}
  
Standard separation of variables implies that the eigenfunctions of $-\Delta^\Dir_{A_r}$, written in polar coordinates $\rho=\sqrt{x_1^2+x_2^2}$ and $\psi$, are of the form 
\[
\left(c_1 J_{m}\left(\sqrt{\lambda}\rho\right)+c_2 Y_{m}\left(\sqrt{\lambda}\rho\right)\right)\er^{\pm\ir m \psi},\qquad m\in\mathbb{N}_0:=\{0\}\cup\mathbb{N},
\]
where $J_m$ and $Y_m$ are the Bessel functions of the first and second kind, respectively, and $|c_1|^2+|c_2|^2>0$. If we define the functions
\begin{equation}\label{eq:Z}
L_{r, m}(x):=J_{m}\left(x\right)Y_{m}\left(r x\right)-Y_{m}\left(x\right)J_{m}\left(r x\right),\qquad m\in\mathbb{N}_0,
\end{equation}
and denote the $k$th positive root of $L_{r, m}(x)$ by $\ell_{r,m,k}$, then the corresponding eigenvalues equal $\ell_{r,m,k}^2$, taken with multiplicity \[
\kappa_m:=\begin{cases}1\qquad&\text{if }m=0.\\
2\qquad&\text{if }m>0.
\end{cases}
\] 

We recall that the Bessel \emph{modulus} function $M_\nu(z)$ and \emph{phase} function $\theta_\nu(x)$ are defined, for $\nu\ge 0$, via
\begin{equation}\label{eq:Mtheta}
J_\nu(x)+\ir Y_\nu(x)=M_\nu(x) \er^{\ir\theta_\nu(x)},\qquad \lim_{x\to 0^+} \theta_\nu(x)=-\frac{\pi}{2}.
\end{equation}
We refer to \cite[\S10.18]{DLMF}, \cite{Hor}, and \cite{FLPSbessel} for more information on the modulus and phase functions; for the moment we only need to recall that 
\[
M_\nu(x)>0\quad\text{and}\quad\theta_\nu'(x)>0\qquad\text{for all }x>0,
\]
and that 
\begin{equation}\label{eq:asympttheta}
\theta_\nu(x) = x - \left(\frac{\nu}{2}+\frac{1}{4}\right)\pi + O\left(\frac{1}{x}\right)\qquad\text{as }x\to+\infty
\end{equation}
\cite[\S10.18.18]{DLMF}.
Using notation \eqref{eq:Mtheta} and definition \eqref{eq:Z}, we get
\[
\begin{split}
L_{r, m}(x)&=M_m(x) M_m(r x) \left(\cos\theta_m(x)\sin\theta_m(r x)-\cos\theta_m(r x)\sin\theta_m(x)\right)\\
&=-M_m(x) M_m(r x)\sin\left(\theta_m(x)-\theta_m(r x)\right).
\end{split}
\]
It is easy to see that the difference 
\begin{equation}\label{eq:deltadef}
\Theta_{r,m}(x):=\theta_m(x)-\theta_m(r x)
\end{equation}
is monotonically increasing in $x$ for every $m\ge 0$ and every $r\in(0,1)$, with 
\[
\Theta_{r,m}(0)=0\qquad\text{and}\qquad\lim_{x\to\infty}\Theta_{r,m}(x)=\infty,
\]
see Lemma \ref{lem:mono}, and therefore the function $\Theta_{r,m}$ is invertible. Thus, the $k$th positive root $\ell_{r,m,k}$ of $L_{r, m}(x)$ is given by
\[
\ell_{r,m,k}=\Theta_{r,m}^{-1}(\pi k),
\]
and the eigenvalue counting function of an annulus can be rewritten as
\begin{equation}\label{eq:Nviadelta}
\mathcal{N}_r(\lambda)=\sum_{m=0}^\infty \kappa_m \#\{k\in\mathbb{N}: \ell_{r,m,k}\le\lambda\}=\sum_{m=0}^\infty \kappa_m \entire{\frac{1}{\pi}\Theta_{r,m}(\lambda)}=\sum_{m=0}^\infty \kappa_m \entire{\frac{1}{\pi}\left(\theta_m(\lambda)-\theta_m(r\lambda)\right)}.
\end{equation}

We will require the following technical results.

\begin{lemma}\label{lem:mono} 
Let $m\ge 0$, $0<r<1$. Then the function $\Theta_{r,m}(x)$ defined by \eqref{eq:deltadef} is strictly monotone increasing at any point $x>0$, and $\Theta_{r,m}(\lambda)\to\infty$ as $\lambda\to\infty$. Moreover, 
\begin{equation}\label{eq:deltalambdalessthanm}
\Theta_{r,m}(\lambda)<\pi\qquad\text{for any }\lambda\in(0,m].
\end{equation}
\end{lemma}

\begin{proof} We have 
\begin{equation}\label{eq:deltaprime}
\Theta'_{r,m}(x)=\theta'_m(x)-r \theta'_m(r x) =\frac{2}{\pi x}\left(\frac{1}{(M_m(x))^2}-\frac{1}{(M_m(r x))^2}\right),
\end{equation}
where the last equality is by \cite[\S10.18.8]{DLMF}. Nicholson's formula  \cite[\S10.9.30]{DLMF},
\[
(M_m(x))^2 = \frac{8}{\pi^2}\int_0^\infty \cosh(2 m t)K_0(2 x \sinh t)\,\dr t,
\]
in which $K_0(\cdot)$ is the modified Bessel function of the second kind (which is strictly decreasing on $(0,+\infty)$), then ensures that $M_m(x)$ is strictly decreasing in $x$, and therefore the right-hand side of \eqref{eq:deltaprime} is positive. 

Also,
\[
\Theta_{r,m}(\lambda) = (1-r)\lambda +O\left(\frac{1}{\lambda}\right)\qquad\text{as }\lambda\to\infty
\]
 by  \eqref{eq:asympttheta}.
 
We also have, for $0<\lambda\le m$,  
\[
\Theta_{r,m}(\lambda)<\theta_m(\lambda)+\frac{\pi}{2}\le \theta_m(m)+\frac{\pi}{2}<\pi,
\] 
since 
\[
-\frac{\pi}{2}<\theta_m(r m)<\theta_m(m)<\theta_m\left(j_{m,1}\right)=\frac{\pi}{2}
\]
by monotonicity of $\theta_m$ and the classical bound $j_{m,1}>m$.
\end{proof}

Lemma \ref{lem:mono} immediately implies that $\ell_{r,m,1}=\Theta_{r,m}^{-1}(\pi)>m$, and thus
\begin{corollary}\label{cor:finite} For an integer $m>\entire{\lambda}\ge 0$,
\[
\#\{k\in\mathbb{N}: \ell_{r,m,k}\le\lambda\}=0.
\]
\end{corollary}

Therefore, \eqref{eq:Nviadelta} can be rewritten  as
\begin{equation}\label{eq:Nviadelta1}
\mathcal{N}_r(\lambda)=\sum_{m=0}^{\entire{\lambda}} \kappa_m \entire{\frac{1}{\pi}\Theta_{r,m}(\lambda)},
\end{equation}
and the statement of Theorem \ref{thm:main} is equivalent to the bound
\begin{equation}\label{eq:bd}
\mathcal{N}_r(\lambda)=\sum_{m=0}^{\entire{\lambda}} \kappa_m \entire{\frac{1}{\pi}\Theta_{r,m}(\lambda)}=\sum_{m=0}^{\entire{\lambda}} \kappa_m \entire{\frac{1}{\pi}\left(\theta_m(\lambda)-\theta_m(r\lambda)\right)}<\frac{\lambda^2(1-r^2)}{4}
\end{equation}
being valid for all  $0<r<1$ and all $\lambda$.

In the sequel, it will be convenient to interchange the independent variable and the parameter, and to switch to $\mu$ instead of $r$ by setting 
\begin{equation}\label{eq:defngamma}
\gamma_{\lambda,\mu}(z):=\frac{1}{\pi}\Theta_{\frac{\mu}{\lambda},z}(\lambda)=\frac{1}{\pi}\left(\theta_z(\lambda)-\theta_z(\mu)\right),
\end{equation}
noting that \eqref{eq:deltalambdalessthanm} in this notation reads
\begin{equation}\label{eq:gammalargelambda}
\gamma_{\lambda,\mu}(m)<1\qquad\text{for all }m\ge \lambda,
\end{equation}
and thus rewriting \eqref{eq:Nviadelta1} as
\begin{equation}\label{eq:Nviadelta2}
\mathcal{N}_r(\lambda)=\sum_{m=0}^{\entire{\lambda}} \kappa_m \entire{\gamma_{\lambda,\mu}(m)}.
\end{equation}

\begin{remark}\label{rem:gamma}
We note that since $\entire{\gamma_{\lambda,\mu}\left(\entire{\lambda}+1\right)}=0$ by \eqref{eq:gammalargelambda}, 
the right-hand side of \eqref{eq:Nviadelta2} can be interpreted as a multiple of the trapezoidal floor sum:  
\begin{equation}\label{eq:NasT}
\mathcal{N}_r(\lambda)=2\mathbf{T}\left(\gamma_{\lambda,\mu}, 0, \entire{\lambda}+1\right).
\end{equation}
\end{remark}

\section{Bounds on the phase functions difference and reduction to a lattice counting problem}\label{sec:gammabounds}

We proceed to estimating the values $\Theta_{r,m}(\lambda)$ appearing in the left-hand side of \eqref{eq:bd}. Introduce, for $\lambda>0$ and $\mu\ge0$,  the functions
\begin{align}
G_\lambda(z)&:=\begin{cases} 
\frac{1}{\pi}\left(\sqrt{\lambda^2-z^2} -z\arccos\frac{z}{\lambda}\right),\quad &z\in[0,\lambda],\\
0,\quad &z>\lambda,
\end{cases}\label{eq:G}\\
H_\mu(z)&:=\frac{3\mu^2+2z^2}{24\pi(\mu^2-z^2)^{3/2}},\qquad\qquad\quad\ z\in[0,\mu),\label{eq:H}\\
F_\mu(z)&:=\begin{cases} 
\max\left\{G_\mu(z)-H_\mu(z),-\frac{1}{4}\right\},\qquad &z\in[0,\mu),\\
-\frac{1}{4},\quad &z\ge \mu.
\end{cases}\label{eq:F}
\end{align}
Also, for $0\le z\le\mu<\lambda$, define
\begin{equation}\label{eq:Phi}
\Phi_{\lambda,\mu}(z):=G_\lambda(z)-G_\mu(z).
\end{equation}

\begin{theorem}\label{thm:thetabound}
Let $z\ge 0$, $\lambda>0$. Then
\[
F_\lambda(z)-\frac{1}{4}<\frac{1}{\pi}\theta_z(\lambda)<G_\lambda(z)-\frac{1}{4}.
\]
\end{theorem}

\begin{proof} 
For $\lambda>z$, the result is the re-statement of \cite[Theorem 1.4]{FLPSbessel} using notation \eqref{eq:G}--\eqref{eq:F}. For $0<\lambda\le z$, the lower bound is the re-statement of the elementary bound $\theta_z(\lambda)>-\frac{\pi}{2}$,  and the upper bound follows from monotonicity and continuous differentiability  of both $\theta_z(\lambda)$ and $G_\lambda(z)$ in $\lambda$ since
\[
\frac{1}{\pi}\theta_z(\lambda)\le \frac{1}{\pi}\theta_z(z)\le \frac{1}{\pi}\lim_{\epsilon\to 0^+} \theta_z(z+\epsilon)\le G_z(z)-\frac{1}{4}=-\frac{1}{4},
\]
and $\theta'_z(\lambda)>0=G_\lambda'(z)$, see also \cite[formula (3.8)]{FLPS}.
\end{proof}

\begin{figure}[ht]
\centering
\includegraphics{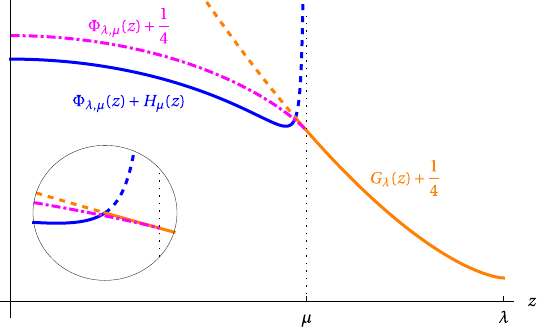}
\caption{A typical behaviour of bounds \eqref{eq:deltabound2} and  \eqref{eq:deltabound3}, and of the upper bound \eqref{eq:deltaboundNF}, with the solid line showing bound \eqref{eq:deltabound1}.  The inset zooms near intersections of the curves.\label{fig:boundsGFH}}
\end{figure}

Theorem \ref{thm:thetabound} plays the key role in the proof of the following bounds, which are illustrated in Figure \ref{fig:boundsGFH}.

\begin{lemma}\label{cor:deltabound}
Let $0<r<1$, $\lambda>0$, $0\le z\le \lambda$, $\mu=r\lambda$, $\gamma_{\lambda,\mu}(z)=\frac{1}{\pi}\Theta_{r, z}(\lambda)$ as before. Then
\begin{equation}\label{eq:deltabound1}
\gamma_{\lambda,\mu}(z)<G_\lambda(z)-F_{\mu}(z).
\end{equation}
In particular,
\begin{equation}\label{eq:deltabound2}
\gamma_{\lambda,\mu}(z)<G_\lambda(z)+\frac{1}{4},
\end{equation}
and if $z<\mu$, then also
\begin{equation}\label{eq:deltabound3}
\gamma_{\lambda,\mu}(z)<\Phi_{\lambda,\mu}(z)+H_\mu(z).
\end{equation}
Additionally, if $z\le\mu$, then also 
\begin{equation}\label{eq:deltaboundNF}
\Phi_{\lambda,\mu}(z)< \gamma_{\lambda,\mu}(z)< \Phi_{\lambda,\mu}(z)+\frac{1}{4}.
\end{equation}
\end{lemma}

\begin{proof}[Proof of Lemma \ref{cor:deltabound}] As $\gamma_{\lambda,\mu}(z)=\frac{1}{\pi}\left(\theta_z(\lambda)-\theta_z(\mu)\right)$ by \eqref{eq:defngamma}, bound \eqref{eq:deltabound1} follows immediately from Theorem \ref{thm:thetabound} by applying its upper bound to $\theta_z(\lambda)$ and its lower bound to $\theta_z(\mu)$. By definitions \eqref{eq:F} and \eqref{eq:Phi}, the right-hand side of \eqref{eq:deltabound1} becomes
\[
G_\lambda(z)-F_{\mu}(z)=
\begin{cases} 
\min\left\{\Phi_{\lambda,\mu}(z)+H_\mu(z), G_\lambda(z)+\frac{1}{4}\right\},\qquad &z\in[0,\mu),\\
G_\lambda(z)+\frac{1}{4},\quad &z\ge \mu,
\end{cases}
\]
thus providing \eqref{eq:deltabound2} and \eqref{eq:deltabound3}.

In order to prove \eqref{eq:deltaboundNF}, fix $z\ge 0$ and consider the function 
\[
\delta_z(\mu):=G_\mu(z)-\frac{1}{\pi}\theta_z(\mu),\qquad \mu\in[z, +\infty).
\]
We have
\[
\delta_z(z)=-\frac{1}{\pi}\theta_z(z)\le -\frac{1}{\pi}\theta_z(0)=\frac{1}{2}
\]
by the equality $G_z(z)=0$ and monotonicity of $\theta$,
\[ 
\lim_{\mu\to+\infty} \delta_z(\mu)=\frac{1}{4}
\]
by \eqref{eq:asympttheta} and the asymptotics
\[
G_\mu(z)=\frac{\mu}{\pi}-\frac{z}{2}+O\left(\frac{1}{\mu}\right)\qquad\text{as }\mu\to+\infty,
\] 
\cite[formula (3.7)]{FLPS}, and also
\[
\delta'_z(\mu) = \frac{1}{\pi}\left(\frac{\sqrt{\mu^2-z^2}}{\mu}-\theta'_z(\mu)\right) <0 \qquad\text{for all }\mu\ge z
\]
by \cite{Hor}. 
Therefore $\frac{1}{4}< \delta_z(\mu)\le \frac{1}{2}$ for all $\mu\in[z,+\infty)$, and thus
\[
0< \gamma_{\lambda,\mu}(z)-\Phi_{\lambda,\mu}(z) = \delta_z(\mu)-\delta_z(\lambda)<\frac{1}{4}
\]
because $\lambda>\mu$.
\end{proof}

In the next section, we summarise some additional properties of the functions \eqref{eq:G}--\eqref{eq:Phi} which we will use afterwards.  

\section{Some properties of functions $G_\lambda$ and $\Phi_{\lambda,\mu}$}\label{sec:aux}

\begin{lemma}\label{lem:G0}
Let $0<z<\lambda$. Then $G_\lambda$ is decreasing and convex on $\left[0, \entire{\lambda}+1\right]$, 
\[
\int_0^{\entire{\lambda}+1}G_\lambda(z)\,\dr z = \int_0^\lambda G_\lambda(z)\,\dr z=\frac{\lambda^2}{8},
\]
and also $G_\lambda$ is $\Lip_{\frac{1}{2}}$ on $\left[0, \entire{\lambda}+1\right]$, and $\Lip_{\frac{1}{3}}$ on  $\left[\frac{\lambda}{2}, \entire{\lambda}+1\right]$, with 
\begin{equation}\label{eq:Glambda2}
G_\lambda\left(\frac{\lambda}{2}\right)=\omega_0\lambda, \qquad \omega_0:=\frac{\sqrt{3}}{2\pi}-\frac{1}{6}\approx 0.108998>\frac{1}{10}.
\end{equation}
\end{lemma}

\begin{proof} The result follows by direct computation and is mostly contained in \cite[Lemmas 4.5 and 4.6]{FLPS}. 
\end{proof}

\begin{lemma}\label{lem:G}
Let $0<z<\lambda$. Then\footnote{We do not actually use the lower bound in this paper.}
\[
\frac{3\lambda}{10}\left(1-\frac{z}{\lambda}\right)^{3/2}<G_\lambda(z)<\frac{\lambda}{3}\left(1-\frac{z}{\lambda}\right)^{3/2}.
\]
\end{lemma}

\begin{proof}
We first prove that for $0<w<\lambda$, 
\begin{equation}\label{eq:wineq}
\frac{9}{20}\sqrt{1-\frac{w}{\lambda}}<-G'_\lambda(w)=\frac{1}{\pi}\arccos\frac{w}{\lambda}<\frac{1}{2}\sqrt{1-\frac{w}{\lambda}},
\end{equation}
or, equivalently,
\begin{equation}\label{eq:vineq}
\cos\left(\frac{\pi}{2}\sqrt{1-v}\right)<v<\cos\left(\frac{9\pi}{20}\sqrt{1-v}\right),
\end{equation}
where $v:=\frac{w}{\lambda}\in(0,1)$. 
The left inequality in \eqref{eq:vineq} holds since $\left.\cos\left(\frac{\pi}{2}\sqrt{1-v}\right)-v\right|_{v\in\{0,1\}}=0$ and 
\[
\frac{\dr^2}{\dr v^2}\left(\cos\left(\frac{\pi}{2}\sqrt{1-v}\right)-v\right) = \frac{\pi\cos\left(\frac{\pi}{2}\sqrt{1-v}\right)}{16(1-v)}\left(\frac{2}{\sqrt{1-v}}\tan\left(\frac{\pi}{2}\sqrt{1-v}\right)-\pi\right)>0
\]
for all $v\in(0,1)$.  The right inequality holds since $\left.v-\cos\left(\frac{9\pi}{20}\sqrt{1-v}\right)\right|_{v=1}=0$ and
\[
\frac{\dr}{\dr v}\left(v-\cos\left(\frac{9\pi}{20}\sqrt{1-v}\right)\right) = 1 - \frac{9\pi\sin\left(\frac{9\pi}{20}\pi\sqrt{1-v}\right)}{40\sqrt{1-v}} > 1-\frac{81\pi^2}{800}>0
\]
for all $v\in(0,1)$.

Integrating \eqref{eq:wineq}  in $w$ from $z$ to $\lambda$ gives the result.
\end{proof}

Lemma \ref{lem:G} implies
\begin{corollary}\label{cor:Gintbound} Let $1\le \mu\le\lambda$. Then
\begin{equation}\label{eq:intGN2bound}
\int_{\entire{\mu}}^{\mu} G_\mu(z)\,\dr z = \int_{\entire{\mu}}^{\mu} G_\lambda(z)\,\dr z -  \int_{\entire{\mu}}^{\mu} \Phi_{\lambda,\mu}(z)\,\dr z<\frac{2}{15\sqrt{\mu}}.
\end{equation}
\end{corollary}

\begin{proof}
The first equality is obvious by the definitions \eqref{eq:G} and \eqref{eq:Phi}. In order to estimate the left-hand side, we just integrate the upper bound in Lemma \ref{lem:G} with $\lambda=\mu$  from $\mu-1<\entire{\mu}$ to $\mu$ arriving at the right-hand side of \eqref{eq:intGN2bound}. 
\end{proof}

\begin{lemma}\label{lem:cmulambda} 
Let $\lambda>\mu>0$. Then $\Phi_{\lambda,\mu}(z)$ is decreasing, concave, and $\Lip_{c_{\lambda,\mu}}$ on $[0, \mu]$, 
with
\begin{equation}\label{eq:cbound}
c_{\lambda,\mu}:=\frac{1}{\pi} \arccos\frac{\mu}{\lambda}<\frac{1}{2}\sqrt{1-\frac{\mu}{\lambda}}<\frac{1}{2}.
\end{equation}
Moreover, the first three  derivatives of $\Phi_{\lambda,\mu}(z)$ are all negative for $z\in(0,\mu)$.
\end{lemma}

\begin{proof} We have
\[
\begin{split}
\Phi'_{\lambda,\mu}(z)&=\frac{1}{\pi}\left(\arcsin\frac{z}{\lambda}-\arcsin\frac{z}{\mu}\right)<0,
\\
\Phi''_{\lambda,\mu}(z)&=\frac{1}{\pi}\left(\frac{1}{\sqrt{\lambda^2-z^2}}-\frac{1}{\sqrt{\mu^2-z^2}}\right)<0,
\\
\Phi'''_{\lambda,\mu}(z)&=\frac{z}{\pi}\left(\left(\lambda^2-z^2\right)^{-3/2}-\left(\mu^2-z^2\right)^{-3/2}\right)<0.
\end{split}
\]
We also have
\[
c_{\lambda,\mu}=-\Phi'_{\lambda,\mu}(\mu)=\frac{1}{2}-\frac{1}{\pi}\arcsin\frac{\mu}{\lambda}=\frac{1}{\pi} \arccos\frac{\mu}{\lambda},
\]
and the upper bound in \eqref{eq:cbound} follows from the upper bound in \eqref{eq:wineq}.
\end{proof}

We will later use the following
\begin{corollary}\label{cor:Phidiff}
Let $\lambda>\mu>0$. Then 
\[
\Phi_{\lambda,\mu}(0)-\Phi_{\lambda,\mu}(z)>\frac{(\lambda-\mu)z^2}{2\pi \lambda \mu}
\]
for all $z\in(0,\mu)$.
\end{corollary}

\begin{proof}
We have 
\[
\Phi_{\lambda,\mu}' (0) = 0, \qquad \Phi''_{\lambda,\mu}(0)=-\frac{\lambda-\mu}{\pi \lambda \mu},
\]
and, by Taylor's theorem with the remainder and by Lemma \ref{lem:cmulambda},
\[
\Phi_{\lambda,\mu}(0)-\Phi_{\lambda,\mu}(z)>\frac{z^2}{2}\left(-\Phi''_{\lambda,\mu}(0)\right)=\frac{(\lambda-\mu)z^2}{2\pi \lambda \mu}
\]
for all $z\in (0,\mu)$.
\end{proof}

\section{Lattice point count and the trapezoidal floor sums}\label{sec:apptrap}

\subsection{The outline of the scheme}
Our approach will be as follows. Recall that by Remark \ref{rem:gamma} we need to prove that
\[
\mathcal{N}_r(\lambda)=2\mathbf{T}\left(\gamma_{\lambda,\mu}, 0, \entire{\lambda}+1\right)<\frac{\lambda^2-\mu^2}{4}.
\]
Depending on the values of $\lambda$ and $\mu$, we will use the upper bounds on $\gamma_{\lambda,\mu}$ obtained in \S\ref{sec:gammabounds} and some further bounds, sometimes applied piecewise, to replace the function $\gamma_{\lambda,\mu}(z)$ by some more explicit function $\tilde{\gamma}_{\lambda,\mu}(z)\ge \gamma_{\lambda,\mu}(z)$. Then we will use the estimates of trapezoidal floor sums from \S\ref{sec:sums} to select the values of $\lambda, \mu$ for which we can guarantee that 
\begin{equation}\label{eq:defnDelta}
\frac{\lambda^2-\mu^2}{4}-2\mathbf{T}\left(\tilde{\gamma}_{\lambda,\mu}, 0, \entire{\lambda}+1\right)>0.
\end{equation}
which in turn immediately implies \eqref{eq:bd}.

We also recall, for further use, that
\[
2\left(\int_0^\mu \Phi_{\lambda,\mu}(z)\,\dr z + \int_\mu^\lambda  G_{\lambda}(z)\,\dr z\right) = 2\left(\int_0^\lambda G_{\lambda}(z)\,\dr z - \int_0^\mu G_{\mu}(z)\,\dr z\right)=\frac{\lambda^2-\mu^2}{4},
\]
therefore \eqref{eq:defnDelta} holds if 
\begin{equation}\label{eq:defnDelta1}
\boldsymbol{\Delta}(\lambda,\mu):=\int_0^\mu \Phi_{\lambda,\mu}(z)\,\dr z + \int_\mu^\lambda  G_{\lambda}(z)\,\dr z-\mathbf{T}\left(\tilde{\gamma}_{\lambda,\mu}, 0, \entire{\lambda}+1\right)
\end{equation}
is positive. We will write $\boldsymbol{\Delta}_\aleph(\lambda,\mu)$ when estimating $\boldsymbol{\Delta}(\lambda,\mu)$ in region $\LM_\aleph$. 

\subsection{Region III}\label{subsec:regionIII}

We use \eqref{eq:deltabound2} to set
\[
\tilde{\gamma}_{\lambda,\mu}(z):=G_\lambda(z)+\frac{1}{4},\qquad z\in\left[0, \entire{\lambda}+1\right],
\] 
then
\[
2 \mathbf{T}\left(\gamma_{\lambda,\mu},0,\entire{\lambda}+1\right) \le 2\mathbf{T}\left(G_\lambda+\frac{1}{4}, 0, \entire{\lambda}+1\right).
\]
Recalling Lemma \ref{lem:G0}, we apply Theorem \ref{thm:conveximproved} with $g(z)=G_\lambda(z)$, $[a,b]=[0,\entire{\lambda}+1]$, to deduce that
\begin{equation}\label{eq:TG}
2\mathbf{T}\left(G_\lambda+\frac{1}{4}, 0, \entire{\lambda}+1\right)<\frac{\lambda^2}{4}-\frac{2\entire{\omega_0\lambda}}{4},
\end{equation}
and \eqref{eq:bd} follows if
\[
2\boldsymbol{\Delta}_\mathrm{III}(\lambda,\mu)=\frac{\lambda^2-\mu^2}{4}-2\mathbf{T}\left(G_\lambda+\frac{1}{4}, 0, \entire{\lambda}+1\right)>\frac{1}{4}\left(2\entire{\omega_0\lambda}-\mu^2\right)\ge 0.
\]
As $\entire{\omega_0\lambda}>\omega_0\lambda-1$, we have, using additionally \eqref{eq:Glambda2},
\[
\boldsymbol{\Delta}_\mathrm{III}(\lambda,\mu)>\frac{1}{8}\left(2\omega_0\lambda-2-\mu^2\right)>\frac{1}{8}\left(\frac{\lambda}{5}-2-\mu^2\right),
\]
and keeping the right-hand side non-negative, we obtain

\begin{theorem}\label{thm:III}
Set
\[
\eta_{\mathrm{III},\pm}(r):=\frac{1\pm \sqrt{1-200 r^2}}{10 r^2},\qquad\zeta_\mathrm{III}(\lambda):=\sqrt{\frac{\lambda}{5}-2}.
\]
Inequality \eqref{eq:PolyaB}  holds for all $(\lambda,\mu)\in\LM_\mathrm{III}$, where
\[
\LM_\mathrm{III}:=\left\{(\lambda,\mu): \lambda>10,\ 0<\mu\le\zeta_\mathrm{III}(\lambda)\right\}\subset\LM.
\]
Equivalently,  inequality \eqref{eq:PolyaA}  holds for  all $(r,\lambda)\in\RL_\mathrm{III}$, where
\begin{equation}\label{eq:reg3r}
\RL_\mathrm{III}:=\left\{(r,\lambda): 0<r< \frac{1}{10\sqrt{2}},\ \eta_{\mathrm{III},-}(r)\le\lambda\le \eta_{\mathrm{III},+}(r)\right\}\subset\RL.
\end{equation}
\end{theorem}

\begin{proof} If  $\lambda>10$ and $\mu\le \sqrt{\frac{\lambda}{5}-2}$, then, obviously, $\boldsymbol{\Delta}_\mathrm{III}(\lambda,\mu)>0$. Rewriting this conditions in terms of $r$ and $\lambda$ yields \eqref{eq:reg3r}.
\end{proof} 

We also immediately obtain the improved version of P\'olya's conjecture for the disk.

\begin{proof}[Proof of Theorem \ref{thm:diskimproved}]
By \cite[Theorem 2.3]{FLPS}, re-written using the notation of the present paper, 
\[
\mathcal{N}^\Dir_{\mathbb{D}}(\lambda)\le 2\mathbf{T}\left(G_\lambda+\frac{1}{4}, 0, \entire{\lambda}+1\right),
\]
and hence inequality \eqref{eq:TG} directly implies \eqref{eq:improvedPolya}.
\end{proof}

\subsection{Sharper bounds: Region IV}\label{sec:regIV} 
Let $N_2:=\entire{\mu}$ and $N_3:=\entire{\lambda}+1$. We have, by \eqref{eq:deltaboundNF} and \eqref{eq:deltabound2},
\[
\gamma_{\lambda,\mu}(z)< \tilde{\gamma}_{\lambda,\mu}(z):=
\begin{cases}
\Phi_{\lambda,\mu}(z)+\frac{1}{4}\quad&\text{if }z\in[0, N_2],\\
G_{\lambda}(z)+\frac{1}{4}\quad&\text{if }z\in[N_2, N_3].
\end{cases}
\]
Applying Proposition \ref{prop:gconcave} with $g=\Phi_{\lambda,\mu}(z)+\frac{1}{4}$, $[a,b]=[0,N_2]$ we obtain
\begin{equation}\label{eq:int1TIV}
\int_{0}^{N_2}\Phi_{\lambda,\mu}(z)\,\dr z - \mathbf{T}\left(\tilde{\gamma}_{\lambda,\mu}, 0, N_2\right) > -\frac{1}{4}N_2=-\frac{1}{4}\entire{\mu}.
\end{equation}
Assume that
\begin{equation}\label{eq:mulambdahalf}
\mu\le\frac{\lambda}{2}. 
\end{equation}
Applying Theorem \ref{thm:conveximproved} with $g(z)=G_\lambda(z)$, $[a,b]=[N_2,N_3]$, together with Lemma \ref{lem:G0}, we obtain
\begin{equation}\label{eq:int2TIV}
\int_{N_2}^{N_3} G_{\lambda}(z)\,\dr z - \mathbf{T}\left(\tilde{\gamma}_{\lambda,\mu}, N_2, N_3\right) >\frac{1}{4}\entire{\omega_0\lambda}.
\end{equation}
Adding \eqref{eq:int1TIV} and \eqref{eq:int2TIV} and using Corollary  \ref{cor:Gintbound}, bound \eqref{eq:Glambda2}, and definition \eqref{eq:defnDelta1}, we get
\begin{equation}\label{eq:Deltareg4}
\boldsymbol{\Delta}_\mathrm{IV}(\lambda,\mu) > \frac{1}{4}\left(\entire{\frac{\lambda}{10}}-\entire{\mu}-\frac{8}{15\sqrt{\mu}}\right),
\end{equation}
and we  choose the conditions for the right-hand side to be non-negative. 

\begin{theorem}\label{thm:IV} 
Set
\begin{equation}\label{eq:eta4}
\eta_{\mathrm{IV}}(r):= \max\left\{\frac{64}{225r}, \frac{10}{1-10r}\right\}=
\begin{cases}
\frac{64}{225r}\quad&\text{if }\ 0<r\le \frac{32}{1445},\\
 \frac{10}{1-10r}\quad&\text{if }\ \frac{32}{1445}<r<\frac{1}{10},
\end{cases}
\end{equation}
\[
\zeta_\mathrm{IV,-}(\lambda):=\frac{64}{225},\qquad \zeta_\mathrm{IV,+}(\lambda):=\frac{\lambda}{10}-1.
\]
Inequality \eqref{eq:PolyaB}  holds for all $(\lambda,\mu)\in\LM_\mathrm{IV}$, where
\[
\LM_\mathrm{IV}:=\left\{(\lambda,\mu): \lambda\ge \frac{578}{45},\  \zeta_\mathrm{IV,-}(\lambda)\le\mu\le\zeta_\mathrm{IV,+}(\lambda)\right\}\subset\LM.
\]
Equivalently,  inequality \eqref{eq:PolyaA}  holds for  all $(r,\lambda)\in\RL_\mathrm{IV}$, where
\[
\RL_\mathrm{IV}:=\left\{(r,\lambda): 0<r<\frac{1}{10},\ \lambda\ge \eta_{\mathrm{IV}}(r)\right\}\subset\RL.
\]
\end{theorem}

\begin{proof} 
For $(\lambda,\mu)\in\LM_\mathrm{IV}$, the condition
\eqref{eq:mulambdahalf} holds, $\frac{\lambda}{10}\ge \mu+1$, and $\frac{8}{15\sqrt{\mu}}\le 1$, hence $\boldsymbol{\Delta}_\mathrm{IV}(\lambda,\mu) >0$. The equivalent statement 
in terms of $r$ and $\lambda$ easily follows from the change of variables.
\end{proof}

\subsection{Even sharper bounds: Region V}
In what follows, given $\lambda>\mu>0$, we intend to choose an $s\in (0,1)$ and an integer number $N_1$ such that 
\begin{equation}\label{eq:N1N2}
0\le N_1 < N_2 = \entire{\mu}<N_3=\entire{\lambda}+1,
\end{equation}
\begin{equation}\label{eq:gammawiths}
\gamma_{\lambda,\mu}(z)\le \Phi_{\lambda,\mu}(z)+s,\qquad z\in[0, N_1],
\end{equation}
and
\begin{equation}\label{eq:N1jump}
\entire{\Phi_{\lambda,\mu}(N_1+1)+\frac{1}{4}}=\entire{\Phi_{\lambda,\mu}(N_1)+\frac{1}{4}}-1.
\end{equation}
Given such a choice, we have a bound 
\[
\gamma_{\lambda,\mu}(z)< \tilde{\gamma}_{\lambda,\mu}(z):=
\begin{cases}
\Phi_{\lambda,\mu}(z)+s\quad&\text{if }z\in[0, N_1],\\
\Phi_{\lambda,\mu}(z)+\frac{1}{4}\quad&\text{if }z\in[N_1, N_2],\\
G_{\lambda}(z)+\frac{1}{4}\quad&\text{if }z\in[N_2, N_3],
\end{cases}
\]
where the estimates from above on intervals $[N_1, N_2]$ and $[N_2, N_3]$ follow from \eqref{eq:deltaboundNF} and \eqref{eq:deltabound2}, respectively. 

\paragraph{The choice of $N_1$ and $s$, and bounds on the interval  $[0, N_1]$.}  Set 
\[
z_*:=\sqrt{\frac{2\pi\lambda\mu}{\lambda-\mu}}.
\]
By Corollary \ref{cor:Phidiff}, $\Phi_{\lambda,\mu}(z_*)<\Phi_{\lambda,\mu}(0)-1$, and therefore
\[
\entire{\Phi_{\lambda,\mu}(z_*)+\frac{1}{4}}<\entire{\Phi_{\lambda,\mu}(0)+\frac{1}{4}}.
\]

Let
\[
N_1:=\max\left\{m\in[0, z_*)\cap\mathbb{Z}:  \entire{\Phi_{\lambda,\mu}(m)+\frac{1}{4}}=\entire{\Phi_{\lambda,\mu}(0)+\frac{1}{4}}\right\},
\]
then \eqref{eq:N1jump} holds.
It may happen that $N_1=0$, and we will treat this case separately.

Assume for the moment that  $z_*<\mu$ and set 
\[
s:=H_\mu(z_*).
\]

\begin{lemma}\label{lem:int1} 
Suppose that 
\begin{equation}\label{eq:int1cond}
\mu>4\pi\qquad\text{and}\qquad \lambda>\frac{\mu^2}{\mu-4\pi},
\end{equation}
and that $N_1$ and $s$ are defined as above.
Then both bounds \eqref{eq:N1N2} and \eqref{eq:gammawiths} are true. Additionally, if $N_1>0$, then 
\begin{equation}\label{eq:int1T}
\int_0^{N_1}\Phi_{\lambda,\mu}(z)\,\dr z - \mathbf{T}\left(\tilde{\gamma}_{\lambda,\mu},0,N_1\right) > -\frac{1}{3\pi}.
\end{equation}
\end{lemma}

\begin{proof}
By the construction of $N_1$, we have $N_1<z_*$. Moreover,
\[
\mu^2-2z_*^2 = \mu^2-\frac{4\pi\lambda\mu}{\lambda-\mu}=\frac{\mu}{\lambda-\mu}\left(\lambda(\mu-4\pi)-\mu^2\right)>0
\]
by \eqref{eq:int1cond}, and therefore $N_1<z_*<\frac{\mu}{\sqrt{2}}$.

Also, $0<\mu-1<\entire{\mu}$, and thus we will have $N_1<\entire{\mu}$ if $\frac{\mu}{\sqrt{2}}\le \mu-1$, which is satisfied due to  the first inequality \eqref{eq:int1cond}.

Finally, assuming $N_1>0$, using bound \eqref{eq:deltabound3} for $\gamma$ and the monotonicity of the function $H_\mu(z)$, we deduce, for $z\in[0, N_1]$,
\[
\begin{split}
\gamma_{\lambda,\mu}(z)&<\Phi_{\lambda,\mu}(z)+H_\mu(z)<\Phi_{\lambda,\mu}(z)+H_\mu(z_*)=\Phi_{\lambda,\mu}(z)+s\\
&=\Phi_{\lambda,\mu}(z) + \frac{3\mu^2+2z_*^2}{24\pi(\mu^2-z_*^2)^{3/2}}\le\Phi_{\lambda,\mu}(z) + \frac{4\mu^2}{24\pi(\mu^2/2)^{3/2}}=\Phi_{\lambda,\mu}(z) +\frac{\sqrt{2}}{3\pi\mu}
\end{split}
\]
(proving along the way \eqref{eq:gammawiths} and $s\le \frac{\sqrt{2}}{3\pi\mu}$). The function $\Phi_{\lambda,\mu}(z) + s$ is concave for $z\in[0,N_1]$. Applying Proposition \ref{prop:gconcave} with $g(z)=\Phi_{\lambda,\mu}(z) + s$, $[a,b]=[0,N_1]$, yields
\[
\mathbf{T}\left(\tilde{\gamma}_{\lambda,\mu},0,N_1\right)\le \int_0^{N_1}(\Phi_{\lambda,\mu}(z)+s)\,\dr z\le  \int_0^{N_1}\Phi_{\lambda,\mu}(z)\,\dr z +\frac{\sqrt{2}}{3\pi\mu}  N_1
< \int_0^{N_1}\Phi_{\lambda,\mu}(z)\,\dr z+\frac{1}{3\pi},
\]
and the result follows.
\end{proof}

\begin{remark} 
If it happens that $N_1=0$, we formally set $\mathbf{T}\left(\tilde{\gamma}_{\lambda,\mu},0,N_1\right):=0$, and then bound \eqref{eq:int1T}  remains valid.
\end{remark}

\paragraph{Bound on the interval $[N_1, N_2]$. } The function $\Phi_{\lambda,\mu}(z)+\frac{1}{4}$ is decreasing, concave, and $\Lip_{c_{\lambda,\mu}}$ with $c_{\lambda,\mu}<\frac{1}{2}\sqrt{1-\frac{\mu}{\lambda}}$ on the interval $[0,\mu]$, see Lemma \ref{lem:cmulambda}. As $[0,\mu]\supset[N_1, N_2]$ under the conditions of Lemma  \ref{lem:int2}, and $
\tilde{\gamma}_{\lambda,\mu}(z)=\Phi_{\lambda,\mu}(z)+\frac{1}{4}$ for $z\in[N_1, N_2]$, we can apply Theorem \ref{thm:t25} to obtain
\begin{equation}\label{eq:int2T}
\begin{split}
\int_{N_1}^{N_2}\Phi_{\lambda,\mu}(z)\,\dr z - \mathbf{T}\left(\tilde{\gamma}_{\lambda,\mu},N_1, N_2\right) &> (N_2-N_1)\left(\frac{1-c_{\lambda,\mu}}{2}-\frac{1}{4}\right)\\
&>\frac{N_2-N_1}{4}\left(1-\sqrt{1-\frac{\mu}{\lambda}}\right)>0.
\end{split}
\end{equation}
We therefore have

\begin{lemma}\label{lem:int2} 
Assume conditions \eqref{eq:int1cond}. Then  bound \eqref{eq:int2T} holds. If, additionally, $\mu\ge\frac{\lambda}{2}$, then
\begin{equation}\label{eq:int2T1}
\int_{N_1}^{N_2}\Phi_{\lambda,\mu}(z)\,\dr z - \mathbf{T}\left(\tilde{\gamma}_{\lambda,\mu},N_1, N_2\right) > \frac{3-2\sqrt{2}}{8}\mu-\frac{2-\sqrt{2}}{8}.
\end{equation}
\end{lemma}

\begin{proof} We have, as shown in the proof of Lemma \ref{lem:int1}, $N_2-N_1>\mu-1-\frac{\mu}{\sqrt{2}}$. Also, if $\mu\ge\frac{\lambda}{2}$, then $1-\sqrt{1-\frac{\mu}{\lambda}}\ge 1-\frac{1}{\sqrt{2}}$, 
and \eqref{eq:int2T1} now follows from \eqref{eq:int2T} after minor simplifications.
\end{proof}

\paragraph{Bounds on the interval $[N_2, N_3]$. } We consider the two cases $\mu<\frac{\lambda}{2}$ and $\mu\ge \frac{\lambda}{2}$ separately. 

In the former case, we act as in \S\ref{sec:regIV} using Theorem \ref{thm:conveximproved} to get
\[
\int_{N_2}^{N_3} G_{\lambda}(z)\,\dr z - \mathbf{T}\left(\tilde{\gamma}_{\lambda,\mu}, N_2, N_3\right) >\frac{1}{4}\entire{\omega_0\lambda},
\]
which together with Corollary  \ref{cor:Gintbound} and bound \eqref{eq:Glambda2} yields
\begin{equation}\label{eq:int3TV1}
\int_{N_2}^\mu \Phi_{\lambda,\mu}(z)\,\dr z +\int_{\mu}^{N_3} G_{\lambda}(z)\,\dr z - \mathbf{T}\left(\tilde{\gamma}_{\lambda,\mu}, N_2, N_3\right)>\frac{1}{4}\left(\entire{\frac{\lambda}{10}}-\frac{8}{15\sqrt{\mu}}\right)
\end{equation}

In the latter case, we use instead Theorem \ref{thm:convexold} to obtain
\begin{equation}\label{eq:int3TV2}
\int_{N_2}^\mu \Phi_{\lambda,\mu}(z)\,\dr z +\int_{\mu}^{N_3} G_{\lambda}(z)\,\dr z - \mathbf{T}\left(\tilde{\gamma}_{\lambda,\mu}, N_2, N_3\right)>-\frac{2}{15\sqrt{\mu}}.
\end{equation}

\paragraph{Putting everything together.} We now combine the bounds we deduced, subject to appropriate conditions, on intervals $[0, N_1]$, $[N_1, N_2]$, and $[N_2, N_3]$. First of all, we remark that in terms of parameters $(r,\lambda)$ conditions \eqref{eq:int1cond} translate as
\begin{equation}\label{eq:lambdarV}
\lambda>\frac{4\pi}{r(1-r)},
\end{equation}
and this in turn implies $\lambda>16\pi$. Also, the first condition \eqref{eq:int1cond} implies $\frac{2}{15\sqrt{\mu}}<\frac{1}{15\sqrt{\pi}}$.

\begin{theorem}\label{thm:V}
Let 
\begin{equation}\label{eq:etazeta5}
\eta_\mathrm{V}(r):= \frac{4\pi}{r(1-r)},\qquad \zeta_{\mathrm{V},\pm}(\lambda):=\frac{1}{2}\left(\lambda \pm \sqrt{\lambda(\lambda-16 \pi)}\right).
\end{equation}
Inequality \eqref{eq:PolyaA}  holds for all $(r,\lambda)\in\RL_\mathrm{V}$, where
\[
\RL_\mathrm{V}:=\left\{(r,\lambda): 0<r<1,\ \lambda> \eta_\mathrm{V}(r)\right\}\subset\RL.
\]
Equivalently,  inequality \eqref{eq:PolyaB} holds for all  $(\lambda,\mu)\in\LM_\mathrm{V}$, where
\[
\LM_\mathrm{V}:=\left\{(\lambda,\mu): \lambda>16\pi,\  \zeta_{\mathrm{V},-}(\lambda)<\mu< \zeta_{\mathrm{V},+}(\lambda)\right\}\subset\LM.
\]
\end{theorem}

\begin{proof} Assume first,  in addition to conditions  \eqref{eq:int1cond},  that $r=\frac{\mu}{\lambda}<\frac{1}{2}$. Adding together bounds \eqref{eq:int1T}, \eqref{eq:int2T}, and \eqref{eq:int3TV1},
and recalling \eqref{eq:defnDelta1}, we obtain
\[
\boldsymbol{\Delta}_\mathrm{V}(\lambda,\mu) > -\frac{1}{3\pi}+\frac{1}{4}\entire{\frac{\lambda}{10}}-\frac{2}{15\sqrt{\mu}}>\frac{\lambda}{40}-\left(\frac{1}{3\pi}+\frac{1}{4}+\frac{1}{15\sqrt{\pi}}\right),
\]
which is positive whenever $\lambda>40\left(\frac{1}{3\pi}+\frac{1}{4}+\frac{1}{15\sqrt{\pi}}\right)\approx 15.7486$. As mentioned, this is automatically true by \eqref{eq:lambdarV}.

Assume now that $r=\frac{\mu}{\lambda}\ge \frac{1}{2}$. Adding together bounds \eqref{eq:int1T}, \eqref{eq:int2T1}, and \eqref{eq:int3TV2},
and recalling once more \eqref{eq:defnDelta1}, we obtain
\[
\boldsymbol{\Delta}_\mathrm{V}(\lambda,\mu) > -\frac{1}{3\pi}+\frac{3-2\sqrt{2}}{8}\mu-\frac{2-\sqrt{2}}{8}-\frac{2}{15\sqrt{\mu}}>\frac{3-2\sqrt{2}}{8}\mu-\left(\frac{1}{3\pi}+\frac{2-\sqrt{2}}{8}+\frac{1}{15\sqrt{\pi}}\right),
\]
which is positive whenever $\mu>\frac{8}{3-2\sqrt{2}}\left(\frac{1}{3\pi}+\frac{2-\sqrt{2}}{8}+\frac{1}{15\sqrt{\pi}}\right)\approx 10.1153$. As we already imposed the restriction $\mu>4\pi$, this is automatically satisfied.
\end{proof}

\section{Filling the void: a computer-assisted algorithm}
\label{sec:compass}

It remains to prove \eqref{eq:PolyaB} in the region 
\begin{equation}\label{eq:LMdiff}
\LM\setminus\LM_\mathrm{theory}.
\end{equation}

\begin{figure}[ht]
\centering
\includegraphics{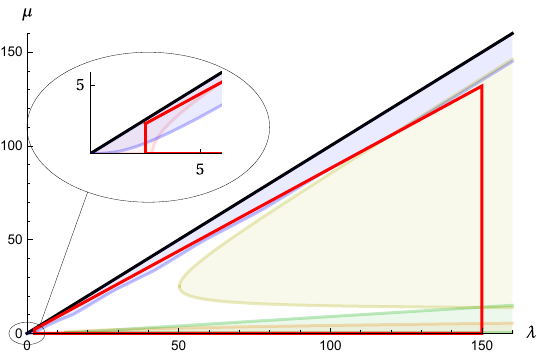}
\caption{The region $\LM_\mathrm{comp}$. \label{fig:LMcomp}}
\end{figure}

\begin{theorem}\label{thm:finiteQ}
The region \eqref{eq:LMdiff} is bounded  and is contained in 
\begin{equation}\label{eq:LMcomp}
\LM_\mathrm{comp} :=\left\{(\lambda,\mu): \frac{5}{2}\le \lambda\le 150, 0\le \mu\le \zeta_\mathrm{comp}(\lambda)\right\},
\end{equation}
where 
\begin{equation}\label{eq:zetacomp}
\zeta_\mathrm{comp}(\lambda):=\frac{22}{25}\lambda,
\end{equation}
see Figure \ref{fig:LMcomp}.
\end{theorem}

\begin{figure}[ht]
\centering
\includegraphics[width=0.45\textwidth]{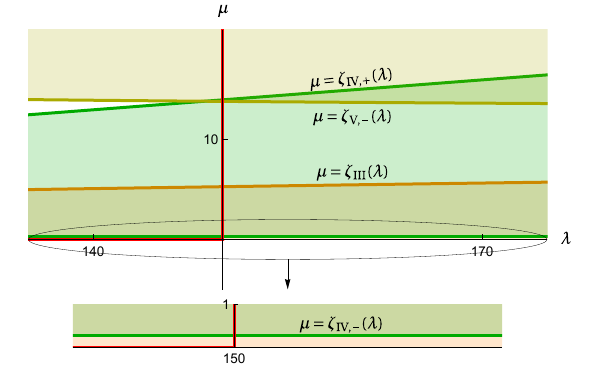}\qquad \includegraphics[width=0.45\textwidth]{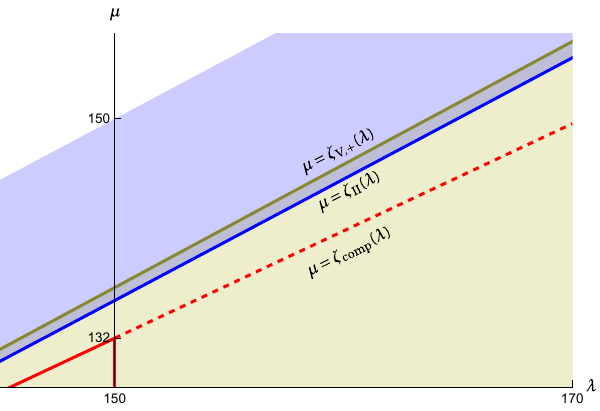}
\caption{The order of the borders of subregions for $\lambda>150$, bottom part on the left (with an extra zoom inset near $\mu=0$),  and top part on the right. \label{fig:thm81c}}
\end{figure}

\begin{proof} We will show that 
\[
(\lambda,\mu)\not\in \LM_\mathrm{comp}\quad\implies\quad(\lambda,\mu)\in\LM_\mathrm{theory}
\]
by considering several cases.

\paragraph{Case 1: $\lambda<\frac{5}{2}$. }  We immediately get $(\lambda,\mu)\in\LM_\mathrm{I}$ by Theorem \ref{thm:regI}.

\paragraph{Case 2: $\mu>\zeta_\mathrm{comp}(\lambda)=\frac{22}{25}\lambda$. } This is the same as $r>\frac{22}{25}$. By \eqref{eq:etazeta5} and \eqref{eq:eta2} (with $j=4$), we have in this case
\[
\frac{\eta_\mathrm{II}(r)}{\eta_\mathrm{V}(r)}=\frac{5 r^{3/2}}{4}>\frac{5 \left(\frac{22}{25}\right)^{3/2}}{4}=\sqrt{\frac{1331}{1250}}>1,
\]
and Theorems \ref{thm:vari} and  \ref{thm:V} imply that either $(r,\lambda)\in\RL_\mathrm{II}$ or $(r,\lambda)\in\RL_\mathrm{V}$, or, equivalently, $(\lambda, \mu)\in\LM_\mathrm{II}\cup\LM_\mathrm{V}$.

\paragraph{Case 3: $\lambda>150$, $0<\mu\le \zeta_\mathrm{comp}(\lambda)=\frac{22}{25}\lambda$ .} We analyse this case by vertical ordering of the boundaries of individual regions in $\LM_\mathrm{theory}$.
We will proceed from bottom to top. Firstly, clearly,
\[
\frac{\zeta_\mathrm{III}(\lambda)}{\zeta_\mathrm{IV,-}(\lambda)}=\frac{225}{64\sqrt{5}}\sqrt{\lambda-10}>\frac{225\sqrt{7}}{32}>1,\qquad 
\frac{\zeta_\mathrm{IV,+}(\lambda)}{\zeta_\mathrm{III}(\lambda)}=\frac{1}{2\sqrt{5}}\sqrt{\lambda-10}>\sqrt{7}>1.
\]
Secondly, by \eqref{eq:eta4} and \eqref{eq:etazeta5},
\[
\frac{\zeta_\mathrm{IV,+}(\lambda)}{\zeta_\mathrm{V,-}(\lambda)}=\frac{\left(1 -\frac{10}{\lambda}\right) \left(5 \lambda +5 \sqrt{\lambda  (\lambda -16 \pi )}\right)}{400\pi},
\]
and since the right-hand side is monotone increasing in $\lambda$, we  have
\[
\frac{\zeta_\mathrm{IV,+}(\lambda)}{\zeta_\mathrm{V,-}(\lambda)}> \frac{\zeta_\mathrm{IV,+}(150)}{\zeta_\mathrm{V,-}(150)}=\frac{7 \left(\sqrt{3 (75-8 \pi )}+15\right)}{60 \pi }\approx 1.01126 > 1.
\]
Thus, by Theorems  \ref{thm:III}, \ref{thm:IV},  and  \ref{thm:V}, $(\lambda,\mu)\in\LM_\mathrm{theory}$ whenever $\lambda>150$ and $0<\mu<\zeta_\mathrm{V,+}(\lambda)$, see the left image in Figure \ref{fig:thm81c}.

Finally, for $\lambda>150$, 
\[
\frac{\zeta_{\mathrm{V},+}(\lambda)}{\zeta_\mathrm{comp}(\lambda)}=\frac{25}{44}\left(1+\sqrt{1-\frac{16}{\pi\lambda}}\right)>\frac{25}{44}\left(1+\sqrt{1-\frac{8}{75\pi}}\right)\approx 1.03148>1,
\]
see the right image in Figure \ref{fig:thm81c}, and therefore   $0<\mu\le \zeta_\mathrm{comp}(\lambda)$ implies $0<\mu<\zeta_\mathrm{V,+}(\lambda)$, completing the proof in Case 3. 
\end{proof}

Let now
\[
\mathcal{P}(\lambda,\mu):=2\mathbf{T}\left(G_\lambda-F_\mu, 0, \entire{\lambda}+1\right),
\]
see \eqref{eq:G}--\eqref{eq:F}. By \eqref{eq:deltabound1} and \eqref{eq:NasT},
\[
\mathcal{N}_{\mu/\lambda}(\lambda)\le \mathcal{P}(\lambda,\mu),
\]
and in order to prove \eqref{eq:PolyaB} for a given $(\lambda,\mu)\in\LM$ it is enough to show that the inequality
\begin{equation}\label{eq:latticeP}
\mathcal{P}(\lambda,\mu)<\frac{\lambda^2-\mu^2}{4}
\end{equation}
holds.

\begin{catheorem}\label{thm:ca} 
Let $(\lambda,\mu)\in  \LM_\mathrm{comp}$. Then \eqref{eq:latticeP} holds.
\end{catheorem}

\begin{remark}\label{rem:cabig} 
We note that the region $\LM_\mathrm{comp}$ is much larger than the region \eqref{eq:LMdiff} which we actually need to deal with, however its shape is determined by the methods we employ. There is, as a result, an element of  redundancy in our computer-assisted algorithm, which covers some parts of the regions already covered by analytic proofs. 
\end{remark}

To start  describing the algorithm of proving Theorem \ref{thm:ca}, we note that $\mathcal{P}(\lambda,\mu)$ is non-decreasing in $\lambda$ and non-increasing in $\mu$: this follows from definitions \eqref{eq:G}--\eqref{eq:F} together with
\begin{align*}
\frac{\partial G_\lambda(z)}{\partial\lambda}&= \frac{\sqrt{\lambda^2-z^2}}{\pi\lambda}>0 \qquad\text{for }0<z<\lambda,\\
\frac{\partial H_\mu(z)}{\partial\mu}&= -\frac{\mu\left(\mu^2+4z^2\right)}{8\pi\left(\mu^2-z^2\right)^{5/2}}<0\qquad\text{for }0<z<\mu.
\end{align*}

To make the algorithm work in a finite number of  exact calculations we use the following two simple but important observations.

\begin{figure}[ht]
\centering
\includegraphics{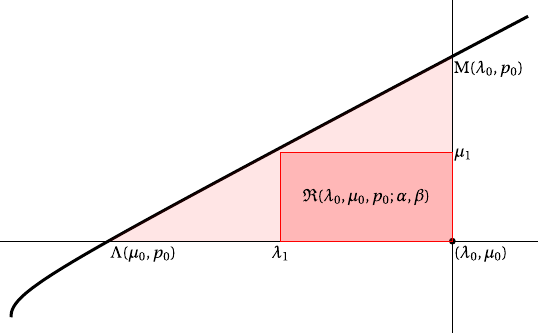}
\caption{The regions \eqref{eq:hyperb} and \eqref{eq:rect}.\label{fig:hyperb}}
\end{figure}

\begin{lemma}\label{lem:rect} Let $0<\mu_0<\lambda_0$, and let $p_0:=\mathcal{P}(\lambda_0,\mu_0)$.
\begin{enumerate}
\item[\normalfont(i)] If $p_0=0$, then \eqref{eq:latticeP} holds in
\[
\mathfrak{T}(\lambda_0,\mu_0):=
\left\{(\lambda,\mu): 0<\lambda\le \lambda_0,  \mu_0\le \mu<\lambda\right\}. 
\]
\item[\normalfont(ii)] If \eqref{eq:latticeP}  holds for $\lambda=\lambda_0$, $\mu=\mu_0$ with a positive \emph{margin}, that is,
\[
\frac{\lambda_0^2-\mu_0^2}{4}-p_0>0,
\]
then it also holds in the following part of the interior of the hyperbola, 
\begin{equation}\label{eq:hyperb}
\left\{(\lambda,\mu): \Lambda(\mu_0, p_0):=\sqrt{\mu_0^2+4p_0}<\lambda\le\lambda_0,\mu_0\le \mu<\sqrt{\lambda^2-4p_0}=:\mathrm{M}(\lambda, p_0)\right\}.
\end{equation}
In particular, it holds inside any rectangle 
\begin{equation}\label{eq:rect}
\mathfrak{R}(\lambda_0, \mu_0, p_0; \alpha, \beta):=[\lambda_1,\lambda_0]\times[\mu_0,\mu_1],
\end{equation}
where $\alpha,\beta\in(0,1)$ and
\[
\begin{split}
\lambda_1&=\lambda_1(\lambda_0, \mu_0, p_0; \alpha):=\alpha\Lambda(\mu_0, p_0)+(1-\alpha)\lambda_0,\\
\mu_1&=\mu_1(\lambda_0, \mu_0, p_0; \alpha,\beta):=\beta\mathrm{M}(\lambda_1, p_0)+(1-\beta)\mu_0,
\end{split}
\]
see Figure \ref{fig:hyperb}.
\end{enumerate}
\end{lemma}

\begin{proof} Statement (i) of the Lemma follows immediately from the monotonicity properties of $\mathcal{P}(\lambda,\mu)$. Again by monotonicity, we
have, for $\lambda\le\lambda_0$ and $\mu\ge \mu_0$,
\[
\frac{\lambda^2-\mu^2}{4}-\mathcal{P}(\lambda,\mu)\ge \frac{\lambda^2-\mu^2}{4}-p_0,
\] 
and this is  positive if $\lambda^2-\mu^2>4p_0$, that is, inside \eqref{eq:hyperb}. It is then elementary to check that  any rectangle $R(\lambda_0, \mu_0, p_0; \alpha, \beta)$ with the restrictions as stated lies inside  \eqref{eq:hyperb}. 
\end{proof}

Lemma \ref{lem:rect} opens a way of covering the region $\LM_\mathrm{comp}$ by a finite number of polygonal regions, always moving up and to the left: these directions explain the reasons for choosing the larger region  for the computer-assisted proof than strictly necessary, see Remark \ref{rem:cabig}.  

Before  describing the algorithm in detail, we also need to address the issue of avoiding floating-point calculations. To do so, we use the method of verified rational approximations we have introduced in \cite{FLPS}. For $x\in\mathbb{R}$, we denote by $ \underline{x}\le x\le \overline{x}$, $\underline{x}, \overline{x}\in\mathbb{Q}$, some rational lower and upper bounds for $x$; see \cite[\S 8]{FLPS} for description of algorithms of finding upper and lower rational approximations of some relevant functions including $G_\lambda$.  Assuming that $\lambda, \mu\in\mathbb{Q}$, we replace the function $\mathcal{P}(\lambda,\mu)$ by the rational-valued function of rational arguments
\[
\overline{\mathcal{P}}(\lambda,\mu)=2\mathbf{T}\left(\overline{G_\lambda}-\underline{F_\mu}, 0, \entire{\lambda}+1\right)\ge \mathcal{P}(\lambda,\mu).
\]

We will construct the finite sequence of rational points, $\left(\lambda^{(k)}, \mu^{(k,j)}\right)$, $k=0,\dots,K$, $j=0,\dots,X_k$, where the number of steps  $K$ and $X_k$ are determined by the algorithm, and the finite sequence of rectangles
\begin{equation}\label{eq:Rkj}
\mathfrak{R}^{(k,j)}:=\left[\lambda^{(k+1)},  \lambda^{(k)}\right]\times\left[\mu^{(k,j)},  \mu^{(k,j+1)}\right],
\end{equation}
(or, occasionally,  triangles), each with its rightmost bottom point at  $\left(\lambda^{(k)}, \mu^{(k,j)}\right)$, inside which \eqref{eq:latticeP} holds, and whose union will cover 
the region $\LM_\mathrm{comp}$.
  
To start with, set $\beta=\frac{99}{100}$ and $\alpha=\frac{2}{3}$ (these values are chosen empirically following some experiments; varying them does not have a radical effect on the required time of computations). Then set 
\[
\LM_\mathrm{comp}^{(0)}:=\LM_\mathrm{comp},\quad 
\zeta_\mathrm{comp}^{(0)}(\lambda):=\zeta_\mathrm{comp}(\lambda),\quad 
\lambda^{(0)} := \max_{(\lambda,\mu)\in\LM_\mathrm{comp}^{(0)}}\lambda= 150,\quad
\mu^{(0,0)} := 0,
\]
and compute 
\[
p^{(0,0)}:= \overline{\mathcal{P}}\left(\lambda^{(0)},\mu^{(0,0)}\right)\in \mathbb{Q},
\]
\[
\lambda^{(1)} := \alpha\overline{\sqrt{{\mu^{(0,0)}}^2+4  p^{(0,0)}}}+(1-\alpha)\lambda^{(0)}\in \mathbb{Q},\quad
\mu^{(0,1)} := \beta\underline{\sqrt{{\lambda^{(1)}}^2-4 p^{(0,0)}}}+(1-\beta)\mu^{(0,0)}\in \mathbb{Q},
\]
and (temporarily) define the rectangle $\mathfrak{R}^{(0,0)}$ by \eqref{eq:Rkj} with $k=j=0$. 
Lemma \ref{lem:rect}, adjusted to rational approximations, then guarantees that  \eqref{eq:latticeP} holds for $(\lambda,\mu)\in \mathfrak{R}^{(0,0)}$.

We now compute a sequence of rectangles $\mathfrak{R}^{(0,j)}$ on top of the first one, by setting on the $j$th vertical step,
\[
p^{(0,j)}:= \overline{\mathcal{P}}\left(\lambda^{(0)},\mu^{(0,j)}\right),\qquad j=1,2,\dots,
\]
\[
\lambda^{(1)}_\mathrm{old}:=\lambda^{(1)},\qquad \lambda^{(1)}_\mathrm{temp}:=\alpha\overline{\sqrt{{\mu^{(0,j)}}^2+4  p^{(0,j)}}}+(1-\alpha)\lambda^{(0)},
\]
\[
\lambda^{(1)}=\lambda^{(1)}_\mathrm{new} := \max\left\{\lambda^{(1)}_\mathrm{temp}, \lambda^{(1)}_\mathrm{old}\right\}
\]
(that is, we always choose the minimal width of all previously and currently computed rectangles in this vertical strip: if $\lambda^{(1)}_\mathrm{new}>\lambda^{(1)}_\mathrm{old}$, we also redefine the previously defined rectangles $\mathfrak{R}^{(0,j')}$, $j'=0,\dots,j-1$, using \eqref{eq:Rkj}), and
\[
\mu^{(0,j+1)} := \beta\underline{\sqrt{{\lambda^{(1)}}^2-4 p^{(0,j)}}}+(1-\beta)\mu^{(0,j)}.
\]
The same argument using Lemma \ref{lem:rect} guarantees that  \eqref{eq:latticeP} holds for $(\lambda,\mu)\in \mathfrak{R}^{(0,j)}$.

\begin{figure}[ht]
\centering
\includegraphics{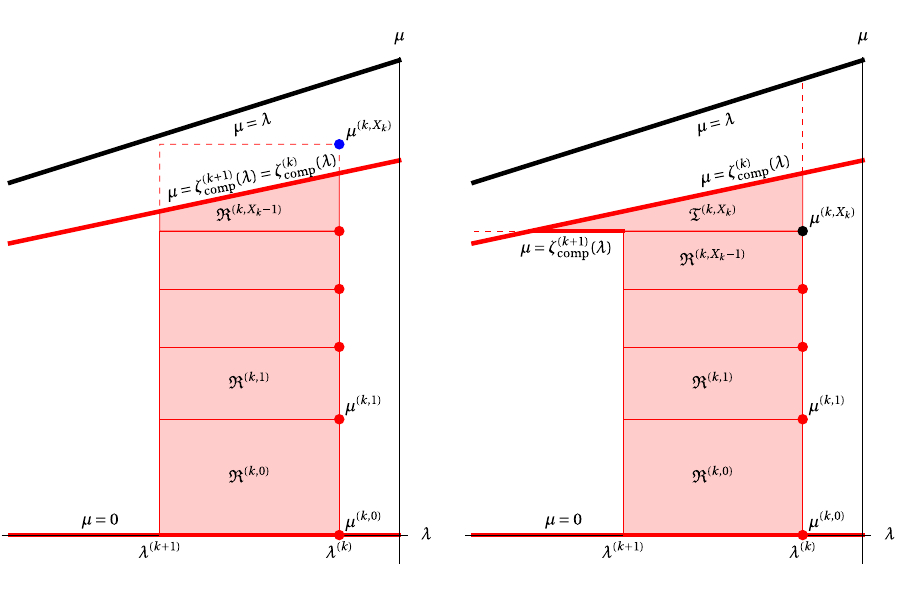}
\caption{The two ways the $k$th  vertical strip may stop: on the left, by reaching a point (blue) outside $\LM_\mathrm{comp}^{(k)}$, and on the right, by reaching a point (black) with the zero trapezoidal floor sum. The image is not to scale. \label{fig:col}}
\end{figure}

The strip is finished in one of the two ways: either for some $j=X_0$ we get outside the region, $\left(\lambda^{(0)}, \mu^{(0,X_0)}\right)\not\in\LM_\mathrm{comp}^{(0)}$, or we get the zero count,  $p^{(0,X_0)}=0$. In both cases we can now remove the vertical strip (or rectangle)
\[
\left\{(\lambda,\mu): \lambda^{(1)}\le \lambda\le \lambda^{(0)}\right\} = \bigcup_{j=0}^{X_0-1} \mathfrak{R}^{(0,j)}
\]
from $\LM_\mathrm{comp}$, and in the latter case we can additionally remove the triangle  $\mathfrak{T}^{\left(0, X_0\right)}=\left\{(\lambda,\mu): \mu\ge\mu^{(0,X_0)}\right\}$, see Figure \ref{fig:col}. In essence we redefine  \eqref{eq:LMcomp}, \eqref{eq:zetacomp} as 
\[
\LM_\mathrm{comp}^{(1)} :=\left\{(\lambda,\mu): \frac{5}{2}\le \lambda\le \lambda^{(1)}, 0\le \mu\le \zeta_\mathrm{comp}^{(1)}(\lambda)\right\},
\]
with
\[
\zeta_\mathrm{comp}^{(1)}(\lambda):=\min\left\{\zeta_\mathrm{comp}^{(0)}(\lambda), \mu^{(0,X_0)}\right\}.
\]
 
We can now proceed to construct another vertical strip, starting from the point $\left(\lambda^{(1)},   \mu^{(1,0)}\right)$,  with $\mu^{(1,0)} := 0$, and so on, with
\[
\LM_\mathrm{comp}^{(k+1)} :=\left\{(\lambda,\mu): \frac{5}{2}\le \lambda\le \lambda^{(k+1)}, 0\le \mu\le \zeta_\mathrm{comp}^{(k+1)}(\lambda)\right\},\quad
\zeta_\mathrm{comp}^{(k+1)}(\lambda):=\min\left\{\zeta_\mathrm{comp}^{(k)}(\lambda), \mu^{(0,X_k)}\right\},
\]
until on some step we reach  $\lambda^{(K)}<\frac{5}{2}$ and have therefore covered all of $\LM_\mathrm{comp}$.

\begin{figure}[ht]
\centering
\includegraphics{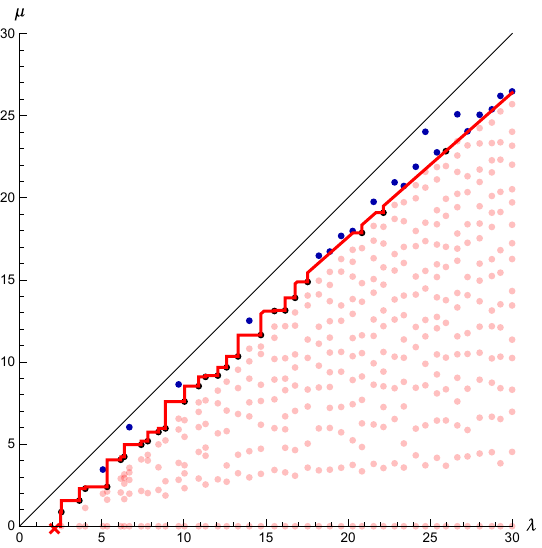}
\caption{The final steps of the calculations. The red line shows the dynamically updated upper boundary of $\LM_\mathrm{comp}^{(k)}$, and the dots are the points where $\overline{\mathcal{P}}$ is evaluated (with blue and black ones as in Figure \ref{fig:col}). The cross is at the point $\left(\lambda^{(K)},0\right)=\left(\frac{179}{82},0\right)$, $K=227$,  where the calculation stops. \label{fig:calcend}}
\end{figure}

The final stages of the calculations are illustrated in Figure \ref{fig:calcend}. Altogether, we complete $K=227$ columns, evaluating the trapezoidal floor sums 8,473 times; the \texttt{Mathematica} script running process required slightly more than 2 minutes of CPU time on a standard Mac laptop.

\begin{figure}[ht]
\centering
\includegraphics{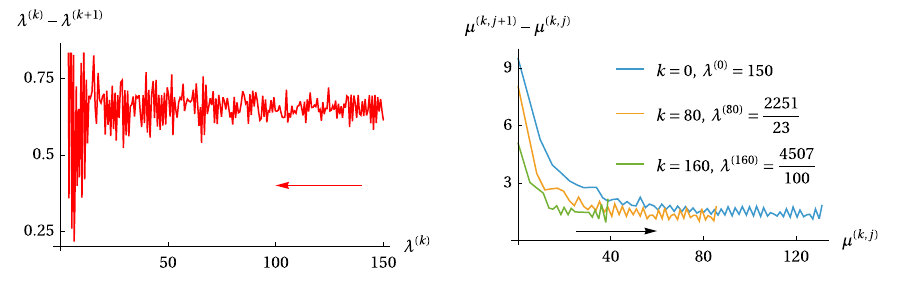}
\caption{On the left, the width of the $k$th vertical strip plotted as a function of its starting abscissa $\lambda^{(k)}$. On the right, the heights of the rectangles $\mathfrak{R}^{(k,j)}$ plotted as functions of $\mu^{(k,j)}$ for some selected $k$s. The arrows indicate the directions of computations.\label{fig:dlambdamu}}
\end{figure}

Some further details on the progress of the computer-assisted algorithm are shown in Figure \ref{fig:dlambdamu}. Full data are available online, see Data availability statement below.

\section*{Data availability statement}\addcontentsline{toc}{section}{Data availability statement}
The accompanying \texttt{Mathematica} script,  its printout, and data files are available for download at \url{https://www.michaellevitin.net/polya.html\#annuli} or at \url{https://github.com/michaellevitin/Polya}.

\section*{Acknowledgements} 
\addcontentsline{toc}{section}{Acknowledgements}
Research of ML was partially supported by EPSRC and by the University of Reading RETF Open Fund. 
Research of IP was partially supported by NSERC and FRQNT.
Research of DS was partially supported by  an AMS-Simons PUI Faculty Research Enhancement Grant.

\phantomsection
\end{document}